\documentclass[a4paper,10pt,twoside]{article}
\usepackage{german,amssymb,amsmath,latexsym,epsfig,graphics,graphicx,mathrsfs,color,amsthm,cite}
\usepackage[latin1]{inputenc}
\usepackage{fancyhdr} 
\newtheorem{theorem}{Theorem}[section]

\newtheorem{lemma}[theorem]{Lemma}

\newtheorem{folgerung}[theorem]{Folgerung}
\newtheorem{beispiel}[theorem]{Beispiel}

\newtheorem{bemerkung}[theorem]{Bemerkung}

\numberwithin{equation}{section}
\newcommand{\R}{{\mathbb R}}
\newcommand{\N}{{\mathbb N}}

\begin{document}
\thispagestyle{empty}
\rule{1.0\textwidth}{1.0pt}

\vspace*{5mm}
{\textbf{\large Optimalit"atsbedingungen f"ur schwache und starke lokale Extrema \\ [2mm]
            in Steuerungsproblemen mit unendlichem Zeithorizont \\[10mm]
            Nico Tauchnitz}} \\[25mm]
{\textbf Vorwort} \\[2mm]
In der vorliegenden Ausarbeitung stelle ich zusammenfassend die Ergebnisse der Arbeiten 
\cite{TauchnitzWMPIHOC,TauchnitzPMPIHOC,TauchnitzOC} zu notwendigen und hinreichenden Optimalit"atsbedingungen f"ur
schwache und starke lokale Extrema in Aufgaben mit unendlichem Zeithorizont vor.
Dabei beginne ich mit Grundlagen zum Dubovickii-Milyutin-Schema,
in dem sich als eine der wesentlichen Kernprobleme die Frage nach der L"osbarkeit linearer Integralgleichungen "uber unbeschr"ankten Intervallen ergibt.
Ich gebe dazu zwei Antworten,
n"amlich einerseits f"ur beschr"ankte Integralkerne im Rahmen gewichteter Funktionenklassen und andererseits f"ur integrable Integralkerne im Rahmen
stetiger Funktionen, die im Unendlichen konvergieren. \\[2mm]
Der Frage nach vern"unftigen L"osungen der linearen Integralgleichung bemesse ich deswegen eine hohe Bedeutung zu,
da die betrachteten Aufgabenklassen mit unendlichem Zeithorizont,
die Anforderungen an die Aufgaben und die erzielten Ergebnisse von dem gew"ahlten Rahmen abh"angen.
Au"serdem gebe ich damit die n"otigen Nacharbeiten der Arbeit \cite{Tauchnitz} an.
In dieser Arbeit sind die linearen Integralgleichungen nicht korrekt behandelt. \\[2mm]
Die vollst"andigen Darstellungen der Ergebnisse sind in \cite{TauchnitzWMPIHOC,TauchnitzPMPIHOC,TauchnitzOC} zu finden.
Insbesondere gehe ich in \cite{TauchnitzWMPIHOC,TauchnitzOC} auf die Schwierigkeiten und Pathologien ein,
die sich im Rahmen gewichteter Funktionenr"aume ergeben k"onnen. \\[5mm]
Juli 2018

\newpage
\lhead[\thepage \, Inhaltsverzeichnis]{Optimale Steuerung mit unendlichem Zeithorizont}
\rhead[Optimale Steuerung mit unendlichem Zeithorizont]{Inhaltsverzeichnis \thepage}
\tableofcontents

\newpage
\lhead[\thepage \hspace*{1mm} Einleitung]{}
\rhead[]{Einleitung \hspace*{1mm} \thepage}

\lhead[\thepage \hspace*{1mm} Dubovickii-Milyutin-Schema]{}
\rhead[]{Dubovickii-Milyutin-Schema \hspace*{1mm} \thepage}
\section{Zur L"osung linearer Integralgleichungen}
Die Herleitung von notwendigen Optimalit"atsbedingungen f"ur das Steuerungsproblem
\begin{equation} \label{EinordnungUH} \left. \begin{array}{l}
  J\big(x(\cdot),u(\cdot)\big) = \displaystyle \int_0^\infty \omega(t)f\big(t,x(t),u(t)\big) \, dt \to \inf, \\[1mm]
  \dot{x}(t) = \varphi\big(t,x(t),u(t)\big), \\[1mm]
  h_0\big(x(0)\big)=0, \qquad \displaystyle\lim_{t \to \infty} h_1\big(t,x(t)\big)=0, \\[1mm]
  u(t) \in U \subseteq \R^m, \quad U \not= \emptyset,
  \end{array} \right\}
\end{equation}
nach dem Dubovickii-Milyutin-Schema (Dubovickii \& Milyutin \cite{DuboMil})
basiert auf der "Uberf"uhrung von (\ref{EinordnungUH}) in eine abstrakte nichtlineare Optimierungsaufgabe der Form
\begin{equation} \label{EA}
J(x,u) \to \inf; \qquad \mathscr{F}(x,u)=0, \qquad x \in X,\; u \in U.
\end{equation}
Dazu betrachten wir an dieser Stelle die Abbildungen
\begin{eqnarray*}
J\big(x(\cdot),u(\cdot)\big) &=& \int_0^\infty \omega(t)f\big(t,x(t),u(t)\big) \, dt, \\
F\big(x(\cdot),u(\cdot)\big)(t) &=& x(t) -x(t_0) -\int_0^t \varphi\big(s,x(s),u(s)\big) \, ds, \quad t \in \R_+,\\
H_0\big(x(\cdot)\big) &=& h_0\big(x(0)\big), \qquad H_1\big(x(\cdot)\big) = \lim_{t \to \infty} h_1\big(t,x(t)\big)
\end{eqnarray*}
und fassen die Gleichrestriktionen zur Abbildung $\mathscr{F}=(F,H_0,H_1)$ zusammen.
Den Kern bei der Herleitung eines Extremalprinzips f"ur die Aufgabe \ref{EA},
d.\,h. beim Nachweis der G"ultigkeit des Lagrangeschen Prinzips,
bildet das Dubovickii-Milyutin-Schema.
Die Vorgehensweise l"asst sich wie folgt abgek"urzt und vereinfacht skizzieren
(die Steuerungsvariable $u(\cdot)$ vernachl"assigen wir):
\begin{enumerate}
\item[(i)] Zu einem lokalen Minimum $x_* \in X$ bezeichne $\mathscr{C}$ die Menge derjenigen $(\eta_0,y)$ mit
           $$\eta_0 > J'(x_*)x, \qquad y = \mathscr{F}'(x_*)x.$$
           Die Menge $\mathscr{C}$ ist konvex und besitzt ein nichtleeres Inneres.
\item[(ii)] W"urde nun der Ursprung dem Inneren der Menge $\mathscr{C}$ angeh"oren,
            so gibt es ein $\overline{x} \in X$ und ferner nach dem Satz von Ljusternik Elemente
            $x_\varepsilon=x_*+ \varepsilon \overline{x}+r(\varepsilon)$ mit
            $$J(x_\varepsilon)<J(x_*), \quad \mathscr{F}(x_\varepsilon)=0, \quad x_\varepsilon \to x_* \mbox{ f"ur } \varepsilon \to 0;$$
            im Widerspruch zur lokalen Optimalit"at von $x_*$.
\item[(iii)] Nach dem Trennungssatz f"ur konvexe Mengen lassen sich die Menge $\mathscr{C}$ und der Ursprung trennen.
             Demnach existiert ein nichttriviales stetiges lineares Funktional (nichttriviale Lagrangesche Multiplikatoren),
             mit dem (denen) das Lagrangesche Prinzip gilt.
\end{enumerate}
Wesentlicher Punkt bei der Anwendung des Trennungssatzes und des Satzes von Ljusternik ist,
dass der Operator $\mathscr{F}'(x_*)$ ein vollst"andiges Bild besitzt
und daher ${\textrm Im\,}F'(x_*)$ mit dem gesamten Bildraum "ubereinstimmt. \\[2mm]
Setzen wir nun zur Abk"urzung $A(t)= \varphi_x\big(t,x(t),u(t)\big)$,
so f"uhrt der eben angegebene Punkt auf die Frage,
unter welchen Rahmenbedingungen die lineare Integralgleichung
$$\zeta(t)=x(t) - \int_0^t A(s)x(s) \, ds, \quad t \in \R_+,$$
zu einer gegebenen Funktion $\zeta(\cdot)$ stets eine L"osung $x(\cdot)$ besitzt.
Eine Antwort auf diese Frage liefert der Nachweis der L"osbarkeit der Fixpunktgleichung $x(\cdot) = T\big(x(\cdot)\big)$, wobei der Operator $T$ durch
$$x(\cdot) \to T\big(x(\cdot)\big), \quad
  T\big(x(\cdot)\big)(t) = \zeta(t) + \int_0^t A(s)x(s) \, ds, \quad t \in \R_+,$$
gegeben ist. Wir geben im Folgenden zwei Varianten an. \\[2mm]
Es bezeichne $C_{\lim}(\R_+,\R^n)$ den Raum der stetigen Funktionen, die im Unendlichen einen Grenzwert besitzen.
Als abgeschlossener Unterraum des Raumes der beschr"ankten stetigen Funktionen ist $C_{\lim}(\R_+,\R^n)$ bez"uglich der Supremumsnorm $\|\cdot\|_\infty$
vollst"andig.

\begin{lemma} \label{FolgerungDGL}
Es sei die Abbildung $t \to A(t)$ "uber $\R_+$ integrierbar.
Dann existiert zu jeder Vektorfunktion $\zeta(\cdot) \in C_{\lim}(\R_+,\R^n)$ eine eindeutig bestimmte Vektorfunktion $x(\cdot) \in C_{\lim}(\R_+,\R^n)$ derart,
dass die lineare Integralgleichung erf"ullt ist:
$$\zeta(t)=x(t) - \int_0^t A(s)x(s) \, ds, \qquad t \in \R_+.$$
\end{lemma}

{\textbf Beweis} Zur abk"urzenden Schreibweise seien 
$$c(t) = \|A(t)\|, \qquad C(t) = \int_0^t c(s) \, ds, \qquad c_0 = \int_0^\infty c(s) \, ds.$$
Bei mehrfacher Anwendung des Operators $T$ ergibt sich f"ur $x_1(\cdot),x_2(\cdot) \in C_{\lim}(\R_+,\R^n)$:
\begin{eqnarray*}
\lefteqn{\big\| \big[T\big(x_1(\cdot) - x_2(\cdot)\big)\big](t) \big\|
         \leq \int_0^t c(s) \, ds \cdot \| x_1(\cdot) - x_2(\cdot) \|_\infty,} \\
\lefteqn{\big\| \big[T^2\big(x_1(\cdot) - x_2(\cdot)\big)\big](t) \big\|
         \leq \int_0^t c(s) \big\| \big[T\big(x_1(\cdot) - x_2(\cdot)\big)\big](s) \big\| \, ds} \\
&& \hspace*{10mm} \leq \int_0^t c(s) C(s) \, ds \cdot \big\| x_1(\cdot) - x_2(\cdot) \big\|_\infty
   = \frac{1}{2} C^2(t) \cdot \| x_1(\cdot) - x_2(\cdot) \|_\infty
\end{eqnarray*}

und sukzessive
\begin{eqnarray*}
\lefteqn{\big\| \big[T^m \big(x_1(\cdot) - x_2(\cdot)\big)\big](t) \big\|
         \leq \int_0^t c(s) \big\| \big[T^{m-1}\big(x_1(\cdot) - x_2(\cdot)\big)\big](s) \big\| \, ds} \\
&& \hspace*{10mm} \leq \int_0^t c(s) \frac{C^{m-1}(t)}{(m-1)!} \, ds \cdot \big\| x_1(\cdot) - x_2(\cdot) \big\|_\infty
       =\frac{C^m(t)}{m!} \cdot \| x_1(\cdot) - x_2(\cdot) \|_\infty.
\end{eqnarray*}
In der Topologie des Raumes $C_{\lim}(\R_+,\R^n)$ gilt daher
$$\big\|T^m \big(x_1(\cdot) - x_2(\cdot)\big) \big\|_\infty \leq \frac{c_0^m}{m!} \cdot \| x_1(\cdot) - x_2(\cdot) \|_\infty.$$
Die Zahlen $a_m = \frac{c_0^m}{m!}$ liefern eine Folge, deren Reihe konvergiert.
Nach dem Fixpunktsatz von Weissinger existiert daher genau ein $x(\cdot)$ mit $x(\cdot) = T x(\cdot)$.
\hfill $\blacksquare$ \\[2mm]
Die Voraussetzung der Integrierbarkeit der Abbildung $t \to A(t)$ in Lemma \ref{FolgerungDGL} ist "au"serst einschr"ankend.
Es entsteht damit die Herausforderungen die Rahmenbedingungen so festzulegen,
dass eine m"oglichst umfassende Klasse an Integralgleichungen behandelt werden kann.
Eine m"ogliche Antwort liefern die gewichteten Funktionenr"aume:
Es sei $\nu(t)=e^{-at}$ mit $a>0$. Dann gilt die im Weiteren wichtige Beziehung
$$\sup_{t \in \R_+} \int_0^t \frac{\nu(t)}{\nu(s)} \, ds = \frac{1}{a}.$$
Mit $C_0(\R_+,\R^n;\nu)$ bezeichnen wir den Raum stetiger Funktionen,
die bez"uglich der Gewichtsfunktion $\nu(\cdot)$ im Unendlichen verschwinden:
$$C_0(\R_+,\R^n;\nu) = \big\{ x(\cdot) \in C(\R_+,\R^n) \,\big|\, \lim_{t \to \infty}\nu(t) x(t) =0\big\}.$$
Wir versehen den Raum $C_0(\R_+,\R^n;\nu)$ mit der gewichteten Supremumsnorm
$$\|x(\cdot)\|_{\infty,\nu} = \sup_{t \in \R_+} \nu(t) \|x(t)\|.$$
Offensichtlich gelten damit die "Aquivalenzen
$$x(\cdot) = \nu^{-1}(\cdot) y(\cdot) \in  C_0(\R_+,\R^n;\nu) \Leftrightarrow \nu(\cdot)x(\cdot) =  y(\cdot) \in C_0(\R_+,\R^n)$$
und $\|x(\cdot)\|_{\infty,\nu} = \|y(\cdot)\|_\infty.$
Damit ist $C_0(\R_+,\R^n;\nu)$ bez"uglich $\|\cdot\|_{\infty,\nu}$ vollst"andig.

\begin{lemma} \label{LemmaDGL2}
Es seien die Abbildung $t \to A(t)$ "uber $\R_+$ me"sbar und beschr"ankt.
Ferner gelte mit der Funktion $\nu(t)=e^{-at}$, $a>0$, die Bedingung
$$\sup_{t \in \R_+} \int_0^t \frac{\nu(t)}{\nu(s)} \|A(s)\| \, ds = q < 1,$$
die z.\,B. f"ur $a > \|A(\cdot)\|_{L_\infty}$ erf"ullt ist.
Dann existiert zu jedem $\zeta(\cdot) \in C_0(\R_+,\R^n;\nu)$ eine eindeutig bestimmte Vektorfunktion $x(\cdot) \in C_0(\R_+,\R^n;\nu)$ derart,
dass die lineare Integralgleichung erf"ullt ist:
$$\zeta(t)=x(t) - \int_0^t A(s)x(s) \, ds, \qquad t \in \R_+.$$
\end{lemma}

{\textbf Beweis} Wir zeigen,
dass der Operator $T$ den Raum $C_0(\R_+,\R^n;\nu)$ in sich abbildet:
Wegen $\zeta(\cdot), x(\cdot) \in C_0(\R_+,\R^n;\nu)$ gelten die Grenzwerte
$$\lim_{t \to \infty} \nu(t) \zeta(t)=0, \qquad \lim_{t \to \infty} \nu(t) x(t)=0.$$
Sei weiterhin $\varepsilon >0$ gegeben.
Dann lassen sich Zahlen $0<T<T'$ angeben mit
$$\nu(t) \|x(t)\| \leq \varepsilon \mbox{ f"ur alle } t \geq T, \qquad \nu(t) \int_0^T \|A(s) x(s)\| \, ds  \leq \varepsilon \mbox{ f"ur alle } t \geq T'.$$
Damit erhalten wir f"ur alle $t \geq T'$
\begin{eqnarray*}
&& \nu(t) \int_0^t \|A(s) x(s)\| \, ds
   = \nu(t) \int_0^T \|A(s) x(s)\| \, ds + \nu(t) \int_T^t \|A(s) x(s)\| \, ds \\
&& \hspace*{10mm} \leq \varepsilon + \|A(\cdot)\|_{L_\infty} \int_T^t \frac{\nu(t)}{\nu(s)} [\nu(s) \|x(s)\|] \, ds
       \leq \varepsilon + \varepsilon \cdot \frac{1}{a} \|A(\cdot)\|_{L_\infty}.
\end{eqnarray*}
Zusammen ergibt sich damit $T\big(x(\cdot)\big)(t) \to 0$ f"ur $t \to \infty$. \\
Wir zeigen nun, dass der Operator $T$ kontraktiv ist:
Bei mehrfacher Anwendung des Operators $T$ ergeben sich f"ur $x_1(\cdot),x_2(\cdot) \in C_0(\R_+,\R^n;\nu)$ die Beziehungen
\begin{eqnarray*}
\lefteqn{\nu(t)\big\| \big[T\big(x_1(\cdot) - x_2(\cdot)\big)\big](t) \big\|
         \leq \nu(t)\int_\tau^t \|A(s)\| \| x_1(s) - x_2(s) \|\, ds} \\
&& \hspace*{10mm} \leq \int_0^t \frac{\nu(t)}{\nu(s)} \big\|A(s)\big\| \cdot \nu(s)\| x_1(s) - x_2(s) \|\, ds
       \leq q \cdot \big\| x_1(\cdot) - x_2(\cdot) \big\|_{\infty,\nu}, \\
\lefteqn{\nu(t)\big\| \big[T^2\big(x_1(\cdot) - x_2(\cdot)\big)\big](t) \big\|
         \leq \nu(t)\int_\tau^t \|A(s)\| \big\| \big[T\big(x_1(\cdot) - x_2(\cdot)\big)\big](s) \big\| \, ds} \\
&& \hspace*{10mm} \leq \int_0^t \frac{\nu(t)}{\nu(s)} \|A(s)\| \cdot \nu(s) \big\|\big[ T\big(x_1(\cdot) - x_2(\cdot)\big)\big](s) \big\| \, ds
       \leq q^2 \cdot \big\| x_1(\cdot) - x_2(\cdot) \big\|_{\infty,\nu}.
\end{eqnarray*}
Sukzessive erhalten wir f"ur $m \in \N$:
\begin{eqnarray*}
\lefteqn{\nu(t) \big\| \big[T^m \big(x_1(\cdot) - x_2(\cdot)\big)\big](t) \big\|
         \leq \nu(t) \int_\tau^t \|A(s)\| \big\| \big[T^{m-1}\big(x_1(\cdot) - x_2(\cdot)\big)\big](s) \big\| \, ds} \\
&& \hspace*{10mm} \leq \int_0^t \frac{\nu(t)}{\nu(s)} \|A(s)\| \cdot \nu(s) \big\| \big[T^{m-1}\big(x_1(\cdot) - x_2(\cdot)\big)\big](s) \big\| \, ds
       \leq q^m \| x_1(\cdot) - x_2(\cdot) \|_{\infty,\nu}.
\end{eqnarray*}
In der Topologie des Raumes $C_0(\R_+,\R^n;\nu)$ gilt daher
$$\big\|T^m \big(x_1(\cdot) - x_2(\cdot)\big) \big\|_{\infty,\nu} \leq q^m \cdot \| x_1(\cdot) - x_2(\cdot) \|_{\infty,\nu}.$$
Wegen $q<1$ existiert nach dem Banachschen Fixpunktsatz eine eindeutige L"osung. \hfill $\blacksquare$

\lhead[\thepage \hspace*{1mm} Schwaches lokales Minimum]{}
\rhead[]{Schwaches lokales Minimum \hspace*{1mm} \thepage}
\section{Schwaches lokales Minimum "uber unendlichem Zeithorizont} \label{KapitelWeak}
\subsection{Die Aufgabenstellung}
Wir untersuchen schwache lokale Minimalstellen der Aufgabe
\begin{eqnarray}
&& \label{WMPAufgabe1} J\big(x(\cdot),u(\cdot)\big) = \int_0^\infty \omega(t) f\big(t,x(t),u(t)\big) \, dt \to \inf, \\
&& \label{WMPAufgabe2} \dot{x}(t) = \varphi\big(t,x(t),u(t)\big), \qquad x(0)=x_0, \\
&& \label{WMPAufgabe3} u(t) \in U \subseteq \R^m, \quad U \not= \emptyset \mbox{ und konvex}, \\
&& \label{WMPAufgabe4} g_j\big(t,x(t)\big) \leq 0 \quad \mbox{f"ur alle } t \in \R_+, \quad j=1,...,l.
\end{eqnarray}
Dabei gelten $f:\R \times \R^n \times \R^m \to \R$, $\varphi:\R \times \R^n \times \R^m \to \R^n$ und $g_j(t,x):\R \times \R^n \to \R$. \\
Wir nennen die Trajektorie $x(\cdot)$ eine L"osung des dynamischen Systems (\ref{WMPAufgabe2}) zur Anfangsbedingung $x(0)=x_0$,
falls $x(\cdot)$ auf $\R_+$ definiert ist und auf jedem endlichen Intervall die Dynamik mit Steuerung $u(\cdot)$
im Sinn von Carath\'eodory l"ost. \\[2mm]
In der Aufgabe (\ref{WMPAufgabe1})--(\ref{WMPAufgabe4}) gilt stets die Annahme:
\begin{enumerate}
\item[(A$_0$)] Es sei $\omega(\cdot) \in L_1(\R_+,\R_+)$ und es sei $\nu(t) = e^{-at}$ ein Gewicht mit $a>0$.
\end{enumerate}
F"ur die Aufgabe (\ref{WMPAufgabe1})--(\ref{WMPAufgabe4}) betrachten wir Variationen im Raum gewichteter stetiger Funktionen,
die im Unendlichen verschwindenn.
Da unter den Eigenschaften der Gewichtsfunktion $\nu(t) = e^{-at}$ die Implikationen
$$x(\cdot) \in W^1_2(\R_+,\R^n;\nu) \quad\Rightarrow\quad \nu(\cdot)x(\cdot) \in W^1_1(\R_+,\R^n) \quad\Rightarrow\quad x(\cdot) \in C_0(\R_+,\R^n;\nu)$$
gelten,
formulieren wir wie Pickenhain \cite{Pickenhain} die Aufgabe (\ref{WMPAufgabe1})--(\ref{WMPAufgabe4}) im Rahmen gewichteter Sobolev-R"aume. \\[2mm]
Wir definieren zu $\big(x(\cdot),u(\cdot)\big) \in W^1_2(\R_+,\R^n;\nu) \times L_\infty(\R_+,U)$ die Menge $U_{\gamma,\nu}$ wie folgt:
$$U_{\gamma,\nu}= \{ (t,x,u) \in \R_+ \times \R^n \times \R^m\,|\, \nu(t)\|x-x(t)\| \leq \gamma, \|u-u(t)\| \leq \gamma\}.$$
Dann geh"oren zur Menge $\mathscr{A}_{\textrm Lip}$ diejenigen $\big(x(\cdot),u(\cdot)\big) \in W^1_2(\R_+,\R^n;\nu) \times L_\infty(\R_+,U)$,
f"ur die es eine Zahl $\gamma>0$ derart gibt, dass auf dem Abschluss der Menge $U_{\gamma,\nu}$ gelten:
\begin{enumerate}
\item[(A$_1$)] Die Abbildungen $f(t,x,u)$, $\varphi(t,x,u)$ und $g_j(t,x)$ sind stetig differenzierbar.
\item[(A$_2$)] F"ur alle $(t,x,u) \in U_{\gamma,\nu}$ ein $L(\cdot) \in L_1(\R_+,\R;\omega)$ und ein $C_0>0$ mit
               \begin{eqnarray*}
               && \big\|\big(f(t,x,u),f_u(t,x,u)\big) \big\| \leq L(t), \quad \|f_x(t,x,u)\| \leq L(t)\nu(t) \\
               && \|\varphi(t,x,u)\| \leq C_0 (1+\|x\|+\|u\|), \quad
                  \big\|\big(\varphi_x(t,x,u),\varphi_u(t,x,u)\big)\big\| \leq C_0.
               \end{eqnarray*}    
               Au"serdem gilt f"ur alle $(t,x,u),(t,x',u') \in U_{\gamma,\nu}$:
               \begin{eqnarray*}
               && \|\varphi_x(t,x',u') - \varphi_x(t,x,u) \| \leq C_0 \big(e^{-at}\|x'-x\| + \|u'-u\|\big), \\
               && \|\varphi_u(t,x',u') - \varphi_u(t,x,u) \| \leq C_0 \big(\|x'-x\| + e^{at} \|u'-u\|\big).
               \end{eqnarray*}
\item[(A$_3$)] F"ur alle $(t,x), (t,x') \in U_{\gamma,\nu}$ existiert ein $C_0>0$ mit
               $$|g_j(t,x)| \leq C_0 (1+\|x\|), \quad \|g_{jx}(t,x)\| \leq C_0, \quad \|g_{jx}(t,x)-g_{jx}(t,x')\| \leq C_0 e^{-at}\|x - x'\|.$$
\end{enumerate}
Die Terme $e^{-at}, e^{at}$ sind eine Konsequenz der Variation in der expandierenden Umgebung $U_{\gamma,\nu}$. \\[2mm]
Wir nennen $\big(x(\cdot),u(\cdot)\big) \in W^1_2(\R_+,\R^n;\nu) \times L_\infty(\R_+,U)$
einen zul"assigen Steuerungsprozess in der Aufgabe (\ref{WMPAufgabe1})--(\ref{WMPAufgabe4}),
falls $\big(x(\cdot),u(\cdot)\big)$ dem System (\ref{WMPAufgabe2}) gen"ugt, die Zustandsbeschr"ankungen (\ref{WMPAufgabe4}) erf"ullt
und das Lebesgue-Integral im Zielfunktional in (\ref{WMPAufgabe1}) endlich ist.
Die Menge $\mathscr{A}_{\textrm adm}$ bezeichnet die Menge der zul"assigen Steuerungsprozesse $\big(x(\cdot),u(\cdot)\big)$. \\[2mm]
Ein zul"assiger Steuerungsprozess $\big(x_*(\cdot),u_*(\cdot)\big)$ ist eine schwache lokale Minimalstelle\index{Minimum, schwaches lokales}
der Aufgabe (\ref{WMPAufgabe1})--(\ref{WMPAufgabe4}),
falls eine Zahl $\varepsilon > 0$ derart existiert, dass die Ungleichung 
$$J\big(x(\cdot),u(\cdot)\big) \geq J\big(x_*(\cdot),u_*(\cdot)\big)$$
f"ur alle $\big(x(\cdot),u(\cdot)\big) \in \mathscr{A}_{\textrm adm}$ mit 
$\|x(\cdot)-x_*(\cdot)\|_\infty \leq \varepsilon$, $\|u(\cdot)-u_*(\cdot)\|_{L_\infty} \leq \varepsilon$ gilt. \\[2mm]
Als abschlie"sende Bemerkung weisen wir darauf hin,
dass die Variationen im gewichteten Rahmen streng genommen auf einen $\nu$-lokalen Optimalit"atsbegriff f"uhren.
D.\,h., dass die Ungleichung
$$J\big(x(\cdot),u(\cdot)\big) \geq J\big(x_*(\cdot),u_*(\cdot)\big)$$
f"ur alle $\big(x(\cdot),u(\cdot)\big) \in \mathscr{A}_{\textrm adm}$ mit 
$\|x(\cdot)-x_*(\cdot)\|_{\infty,\nu} \leq \varepsilon$, $\|u(\cdot)-u_*(\cdot)\|_{L_\infty} \leq \varepsilon$ gilt.
Da aber in der Ungleichung $\|x(\cdot)-x_*(\cdot)\|_{\infty,\nu} \leq \varepsilon$ die Elemente $x(\cdot)$ mit der Eigenschaft
$\|x(\cdot)-x_*(\cdot)\|_\infty \leq \varepsilon$ inbegriffen sind,
entsteht kein Widerspruch.
       
       \subsection{Optimalit"atsbedingungen f"ur eine Grundaufgabe} \label{AbschnittWeak}
Im Weiteren bezeichnet $H: \R \times \R^n \times \R^m \times \R^n \times \R \to \R$ die Pontrjagin-Funktion
$$H(t,x,u,p,\lambda_0) = \langle p, \varphi(t,x,u) \rangle-\lambda_0 \omega(t)f(t,x,u).$$

\begin{theorem} \label{SatzWMP} \index{Schwaches Optimalit"atsprinzip}
Es sei $\big(x_*(\cdot),u_*(\cdot)\big) \in \mathscr{A}_{\textrm adm} \cap \mathscr{A}_{\textrm Lip}$ und $x_*(\cdot) \in W^1_2(\R_+,\R^n;\nu)$.
Ferner sei die folgende Kontraktionsbedingung erf"ullt:
\begin{equation} \label{WMPBedingung}
\sup_{t \in \R_+} \int_0^t \frac{\nu(t)}{\nu(s)} \big\|\varphi_x\big(s,x_*(s),u_*(s)\big)\big\| \, ds < 1.
\end{equation}
Ist $\big(x_*(\cdot),u_*(\cdot)\big)$ ein schwaches lokales Minimum der Aufgabe (\ref{WMPAufgabe1})--(\ref{WMPAufgabe3}),
dann existieren nicht gleichzeitig verschwindende Multiplikatoren $\lambda_0 \geq 0$ und $p(\cdot) \in L_2(\R_+,\R^n;\nu^{-1})$
derart, dass
\begin{enumerate}
\item[(a)] die Funktion $p(\cdot)$ fast "uberall der adjungierten Gleichung\index{adjungierte Gleichung}
           \begin{equation}\label{WMP1}
           \dot{p}(t) = -\varphi_x^T\big(t,x_*(t),u_*(t)\big) p(t) + \lambda_0 \omega(t)f_x\big(t,x_*(t),u_*(t)\big)
           \end{equation}
           gen"ugt und die ``nat"urlichen'' Transversalit"atsbedingungen
           \begin{equation}\label{WMP2}
           \lim_{t \to \infty} \|p(t)\|^2\nu^{-1}(t) = 0, \qquad
           \lim_{t \to \infty} \langle p(t),x(t) \rangle = 0 \quad \forall \; x(\cdot) \in W^1_2(\R_+,\R^n;\nu)
           \end{equation}
           erf"ullt;
\item[(b)] in fast allen Punkten $t \in \R_+$ und f"ur alle $u \in U$ die Variationsungleichung
           \begin{equation}\label{WMP3}
           \big\langle H_u\big(t,x_*(t),u_*(t),p(t),\lambda_0\big),\big(u-u_*(t)\big) \big\rangle\leq 0
           \end{equation}
           gilt.
\end{enumerate}
\end{theorem}

\begin{beispiel}\label{BeispielRegler} {\textrm Wir diskutieren die Voraussetzungen anhand des linear-quadratischen Reglers
$$J\big(x(\cdot),u(\cdot)\big) = \int_0^\infty e^{-\varrho t} \cdot \frac{1}{2}\big( x^2(t)+u^2(t)\big) \, dt \to \inf, \quad
  \dot{x}(t) = \alpha x(t)+\beta u(t).$$
Die Kontraktionsbedingung (\ref{WMPBedingung}) ist nur f"ur ein Gewicht $\nu(t)=e^{-at}$ mit $a > \alpha$ erf"ullt.
D.\,h. ferner, dass die Umgebung
$$U_{\gamma,\nu}= \{ (t,x,u) \in \R_+ \times \R^n \times \R^m\,|\, \nu(t)\|x-x_*(t)\| \leq \gamma, \|u-u(t)\| \leq \gamma\}$$
Trajektorien der Form $x(t) = x_*(t) + \xi(t) e^{at}$ mit $\|\xi(t)\| \leq \gamma$ enth"alt.
Ist dabei die Trajektorie $x_*(\cdot)$ beschr"ankt,
so enth"alt $U_{\gamma,\nu}$ Trajektorien der Form $x(t) = \xi e^{at}$ und es ergibt sich in Voraussetzung (A$_2$)
zwischen der Verteilungs- und Gewichtsfunktion die Bedingung $\varrho > 2a$.
In der Aufgabe
\begin{eqnarray*}
&& J\big(x(\cdot),u(\cdot)\big) = \int_0^\infty e^{-2t} \cdot \frac{1}{2}\big( x^2(t)+u^2(t)\big) \, dt \to \inf, \\
&& \dot{x}(t) = 2 x(t)+u(t), \qquad x(0)=2, \qquad u(t) \in \R
\end{eqnarray*}
liefern die notwendigen Bedingungen (\ref{WMP1})--(\ref{WMP3}) den Steuerungsprozess und die Adjungierte
$$x_*(t)=2e^{(1-\sqrt{2})t}, \quad u_*(t)=-2(1+\sqrt{2})e^{(1-\sqrt{2})t}, \quad p(t)=e^{-2t}u_*(t).$$
Da die Bedingung $\varrho> 2a $ an die Verteilungs- und Gewichtsfunktion verletzt ist,
darf Theorem \ref{SatzWMP} nicht auf den beschr"ankten Steuerungsprozess $\big(x_*(\cdot),u_*(\cdot)\big)$ angewendet werden. \hfill $\square$}
\end{beispiel}

\begin{lemma} \label{LemmaMichel}
Zus"atzlich zu Theorem \ref{SatzWMP} gelte f"ur $L(\cdot) \in L_1(\R_+,\R_+;\omega)$ in (A$_2$):
$$\lim_{t \to \infty} \omega(t)L(t)=0.$$
Dann ist die Bedingung von Michel \cite{Michel}\index{Transversalit"atsbedingungen!von@-- von Michel} erf"ullt:
\begin{equation} \label{NaturalHamilton} \lim_{t \to \infty} H\big(t,x_*(t),u_*(t),p(t),\lambda_0\big)=0. \end{equation}
\end{lemma}

Die folgende ``Stetigkeitsbedingung (S)'' ist hinreichend f"ur die Normalform des Theorems \ref{SatzWMP}:

\begin{enumerate}
\item[(S)] Es existieren ein $T\geq 0$, eine Zahl $\varrho(T)>0$ und ein $\mu_T(\cdot) \in L_2(\R_+,\R_+;\nu)$ derart,
           dass f"ur alle $\zeta_T$ mit $\|\zeta_T-x_*(T)\| \leq \varrho(T)$ das System
           $\dot{x}(t) = \varphi\big(t,x(t),u_*(t)\big)$
           mit Anfangsbedingung $x(T)=\zeta_T$ eine L"osung $x(t;\zeta_T)$ auf $[T,\infty)$ besitzt und folgende Ungleichung gilt
           $$\| x(t;\zeta_T)-x_*(t)\| \leq \|\zeta_T-x_*(T)\|\mu_T(t) \quad \mbox{ f"ur alle } t \geq T \geq 0.$$
\end{enumerate}

\begin{theorem} \label{SatzNormalWMP}
Sei $\big(x_*(\cdot),u_*(\cdot)\big) \in \mathscr{A}_{\textrm adm} \cap \mathscr{A}_{\textrm Lip}$ und seien (\ref{WMPBedingung}), (S) erf"ullt. 
Ist $\big(x_*(\cdot),u_*(\cdot)\big)$ ein schwaches lokales Minimum der Aufgabe (\ref{WMPAufgabe1})--(\ref{WMPAufgabe3}),
dann ist Theorem \ref{SatzWMP} mit $\lambda_0=1$ erf"ullt.
Ferner besitzt die Adjungierte $p(\cdot)$ die Darstellung
\begin{equation} \label{NormalWMP1}
p(t)= -Z_*(t) \int_t^\infty \omega(s) Z^{-1}_*(s) f_x\big(s,x_*(s),u_*(s)\big) \, ds.
\end{equation}
Dabei ist $Z_*(t)$ die in $t=0$ normalisierte Fundamentalmatrix des linearen Systems
$$\dot{z}(t)=-\varphi^T_x\big(t,x_*(t),u_*(t)\big) z(t).$$
\end{theorem}
Die Darstellung (\ref{NormalWMP1}) stimmt
(bis auf das Vorzeichen, das sich durch die Minimierung statt einer Maximierung des Zielfunktionals ergibt)
mit der Formel in den Arbeiten von Aseev \& Kryazhimskii und Aseev \& Veliov \cite{AseKry,AseVel,AseVel2,AseVel3} "uberein.
Im Gegensatz zu diesen Arbeiten ist Theorem \ref{SatzNormalWMP} unter den Voraussetzungen
(A$_0$)--(A$_2$), (\ref{WMPBedingung}), (S) erf"ullt und charakterisiert schwache lokale Minimalstellen.

\begin{theorem} \label{SatzHBWMP}
In der Aufgabe (\ref{WMPAufgabe1})--(\ref{WMPAufgabe3}) sei
$\big(x_*(\cdot),u_*(\cdot)\big)$ ein zul"assiger Steuerungsprozess,
f"ur den die Abbildungen $f(t,x,u)$ und $\varphi(t,x,u)$ auf der Menge
$$U_\gamma= \{ (t,x,u) \in \R_+ \times \R^n \times \R^m\,|\, \|x-x_*(t)\| \leq \gamma, \|u-u_*(t)\| \leq \gamma\}$$
stetig differenzierbar sind.
Ferne gelte:
\begin{enumerate}
\item[(a)] Das Tripel $\big(x_*(\cdot),u_*(\cdot),p(\cdot)\big)$
           erf"ullt (\ref{WMP1})--(\ref{WMP3}) mit $\lambda_0=1$ in Theorem \ref{SatzWMP}.        
\item[(b)] F"ur jedes $t \in \R_+$ ist die Funktion $H\big(t,x,u,p(t),1\big)$
           konkav in $(x,u)$ auf $U_\gamma$.
\end{enumerate}
Dann ist $\big(x_*(\cdot),u_*(\cdot)\big)$ ein schwaches lokales Minimum der Aufgabe (\ref{WMPAufgabe1})--(\ref{WMPAufgabe3}).
\end{theorem}

\begin{beispiel}\label{BeispielDockner} {\textrm Wir betrachten nach Dockner et\,al. \cite{Dockner} das Differentialspiel
\begin{eqnarray*}
&& \tilde{J}_i\big(x(\cdot),u_1(\cdot),u_2(\cdot)\big) =\int_0^\infty e^{-\varrho t}\big(p x(t)-c_i\big)u_i(t) \, dt \to \sup, \\
&& \dot{x}(t)=x(t) \big(\alpha-r\ln x(t) \big) -u_1(t)x(t)-u_2(t)x(t), \quad x(0)=x_0>0, \\
&& u_i > 0, \quad \alpha,c_i,p,r,\varrho >0, \quad \alpha> \frac{1}{c_1+c_2}, \quad i=1,2.
\end{eqnarray*}
Wieder ist der Preis $p$ nicht konstant, sondern umgekehrt proportional zum Angebot:
$$p=p(u_1x+u_2x)= \frac{1}{u_1x+u_2x}.$$
Nach Anwendung der Transformation $z=\ln x$ ergibt sich das Spielproblem
\begin{eqnarray*}
&& J_i\big(z(\cdot),u_1(\cdot),u_2(\cdot)\big)
        =\int_0^\infty e^{-\varrho t}\bigg(\frac{1}{u_1(t)+u_2(t)}-c_i\bigg)u_i(t) \, dt \to \sup, \\
&& \dot{z}(t)=-r z(t)+\alpha-u_1(t)-u_2(t), \quad z(0)=\ln x_0>0, \\
&& u_i>0, \quad \alpha,c_i,p,r,\varrho >0, \quad \alpha> \frac{1}{c_1+c_2}, \quad i=1,2.
\end{eqnarray*}
In diesem Spielproblem sind die Dynamiken linear und au"serdem flie"st die Zustandsvariable nicht im Integranden ein.
W"ahlen wir nun $\nu(t) = e^{-at}$ mit $a >r$,
so sind s"amtliche Annahmen in Theorem \ref{SatzWMP}, insbesondere die Kontraktionsbedingung (\ref{WMPBedingung}), erf"ullt. \\
Der L"osungsansatz "uber ein Nash-Gleichgewicht liefert die Steuerungen
$$u_1^*(t)\equiv \frac{c_2}{(c_1+c_2)^2}, \qquad u_2^*(t)\equiv \frac{c_1}{(c_1+c_2)^2},$$
die optimale Trajektorie
$$z_*(t)=(z_0-c_0)e^{-rt} + c_0, \qquad c_0=  \frac{1}{r}\bigg(\alpha-\frac{1}{c_1+c_2}\bigg)$$
und die Adjungierten
$$p_i(t) \equiv 0, \qquad i=1,2.$$
Die Funktion $z_*(\cdot)$ ist streng monoton und nimmt nur Werte des Segments $[z_0,c_0]$ an.
Da $c_0$ und $z_0$ positiv sind, ist $x_*(t)=\exp\big(z_*(t)\big)$ "uber $\R_+$ wohldefiniert. \hfill $\square$}
\end{beispiel}

\begin{beispiel} \label{BeispielFischereispiel}
{\textrm Im Differentialspiel 
im Beispiel \ref{BeispielDockner} mit unendlichem Zeithorizont "uberf"uhrten wir den Ansatz eines Nash-Gleichgewichtes in die gekoppelten Aufgaben
\begin{eqnarray*}
&& J_i\big(z(\cdot),u_1(\cdot),u_2(\cdot)\big)
        =\int_0^\infty e^{-\varrho t}\bigg(\frac{1}{u_1(t)+u_2(t)}-c_i\bigg)u_i(t) \, dt \to \sup, \\
&& \dot{z}(t)=-r z(t)+\alpha-u_1(t)-u_2(t), \quad z(0)=\ln x_0>0, \\
&& u_i>0, \quad \alpha,c_i,p,r,\varrho >0, \quad \alpha> \frac{1}{c_1+c_2}, \quad i=1,2.
\end{eqnarray*}
Die lineare Dynamik mit Wachstumskoeffizient $-r$ f"uhrt zur Stetigkeitsbedingung
$$\| z(t;\zeta_T)-z_*(t)\| = \|\zeta_T-z_*(T)\|e^{-r(t-T)} = \|\zeta_T-z_*(T)\|\mu_T(t) \quad \mbox{ f"ur alle } t \geq T \geq 0.$$
Dabei geh"ort die Funktion $\mu_T(t)=e^{rT}e^{-rt}$ stets dem Raum $L_2(\R_+,\R_+;\nu)$ mit Gewicht $\nu(t)=e^{-at}$, $a>0$, an.
In beiden Aufgaben ist stets $f_x\big(t,x,u_1,u_2\big) = 0$ und es liefert (\ref{NormalWMP1}) unmittelbar $p_i(t) \equiv 0$ f"ur $i=1,2$. \hfill $\square$}
\end{beispiel}

\begin{beispiel}\label{BeispielRegler2} {\textrm Im linear-quadratischen Regler
\begin{eqnarray*}
&& J\big(x(\cdot),u(\cdot)\big) = \int_0^\infty e^{-2t} \cdot \frac{1}{2}\big( x^2(t)+u^2(t)\big) \, dt \to \inf, \\
&& \dot{x}(t) = 2 x(t)+u(t), \qquad x(0)=2, \qquad u(t) \in \R
\end{eqnarray*}
gen"ugt der Steuerungsprozess und die Adjungierte
$$x_*(t)=2e^{(1-\sqrt{2})t}, \quad u_*(t)=-2(1+\sqrt{2})e^{(1-\sqrt{2})t}, \quad p(t)=e^{-2t}u_*(t)$$
den notwendigen Bedingungen (\ref{WMP1})--(\ref{WMP3}) in Theorem \ref{SatzWMP}.
Weiterhin ist $H\big(t,x,u,p(t),1\big)$ konkav bez"uglich $(x,u)$ auf $U_\gamma$ f"ur jedes $\gamma>0$.
Damit ist $\big(x_*(\cdot),u_*(\cdot)\big)$ ein schwaches lokales Minimum dieser Aufgabe. \hfill $\square$}
\end{beispiel}
       \subsection{Optimalit"atsbedingungen unter Zustandsbeschr"ankungen} \label{AbschnittZustandWMP}
Im Weiteren bezeichne $H: \R \times \R^n \times \R^m \times \R^n \times \R \to \R$ wieder die Pontrjagin-Funktion
$$H(t,x,u,p,\lambda_0) = -\lambda_0 \omega(t)f(t,x,u) + \langle p, \varphi(t,x,u) \rangle.$$

\begin{theorem}\label{SatzWMPAufgabeZB} \index{Schwaches Optimalit"atsprinzip}
Sei $\big(x_*(\cdot),u_*(\cdot)\big) \in \mathscr{A}_{\textrm adm} \cap \mathscr{A}_{\textrm Lip}$.
Weiterhin sei die Kontraktionsbedingung
\begin{equation} \label{WMPBedingungZB}
\sup_{t \in \R_+} \int_0^t \frac{\nu(t)}{\nu(s)} \big\|\varphi_x\big(s,x_*(s),u_*(s)\big)\big\| \, ds < 1
\end{equation}
erf"ullt und es m"ogen f"ur alle $x(\cdot) \in C_0(\R_+,\R^n;\nu)$ mit $\|x(\cdot)-x_*(\cdot)\|_{\infty,\nu} < \gamma$ die folgenden Grenzwerte existieren:
\begin{equation} \label{WMPBedingungZB2}
\lim_{t \to \infty} g_{jx}\big(t,x(t)\big), \qquad j=1,...,l.
\end{equation}
Ist $\big(x_*(\cdot),u_*(\cdot)\big)$ ein schwaches lokales Minimum der Aufgabe (\ref{WMPAufgabe1})--(\ref{WMPAufgabe4}),
dann existieren eine Zahl $\lambda_0 \geq 0$, ein Vektor $l_0 \in \R^n$, eine Funktion $p(\cdot):\R_+ \to \R^n$
und auf den Mengen
$$T_j=\big\{t \in \overline{\R}_+ =\R_+ \cup \{\infty\} \,\big|\, g_j\big(t,x_*(t)\big)=0\big\}, \quad j=1,...,l,$$
konzentrierte nichtnegative regul"are Borelsche Ma"se $\mu_j$ endlicher Totalvariation
(wobei s"amtliche Gr"o"sen nicht gleichzeitig verschwinden) derart, dass
\begin{enumerate}
\item[(a)] die Vektorfunktion $p(\cdot)$ von beschr"ankter Variation ist, der adjungierten Gleichung\index{adjungierte Gleichung}
           \begin{equation}\label{SatzWMPAufgabeZB1}
           p(t)= \int_t^\infty H_x\big(s,x_*(s),u_*(s),p(s),\lambda_0\big) \, ds
                  -\sum_{j=1}^l \int_t^\infty \nu(s) g_{jx}\big(s,x_*(s)\big) \, d\mu_j(s)
           \end{equation}
           gen"ugt und die ``nat"urlichen'' Transversalit"atsbedingungen\index{Transversalit"atsbedingungen!nat@--, nat"urliche}
           \begin{equation}\label{SatzWMPAufgabeZB2}
           \lim_{t \to \infty} \|p(t)\|^2\nu^{-1}(t) = 0, \qquad
           \lim_{t \to \infty} \langle p(t),x(t) \rangle = 0 \quad \forall \; x(\cdot) \in W^1_2(\R_+,\R^n;\nu)
           \end{equation}
           erf"ullt;
\item[(b)] in fast allen Punkten $t \in \R_+$ und f"ur alle $u \in U$ die Variationsungleichung
           \begin{equation}\label{SatzWMPAufgabeZB3}
           \big\langle H_u\big(t,x_*(t),u_*(t),p(t),\lambda_0\big),\big(u-u_*(t)\big) \big\rangle\leq 0
           \end{equation}
           gilt.
\end{enumerate}
\end{theorem}

\begin{theorem} \label{SatzHBWMPZB}
In der Aufgabe (\ref{WMPAufgabe1})--(\ref{WMPAufgabe4}) sei
$\big(x_*(\cdot),u_*(\cdot)\big)$ ein zul"assiger Steuerungsprozess,
f"ur den die Abbildungen $f(t,x,u)$ und $\varphi(t,x,u)$ auf der Menge
$$U_\gamma= \{ (t,x,u) \in \R_+ \times \R^n \times \R^m\,|\, \|x-x_*(t)\| \leq \gamma, \|u-u_*(t)\| \leq \gamma\}$$
stetig differenzierbar sind.
Au"serdem sei die Vektorfunktion $p(\cdot) \in L_2(\R_+,\R^n;\nu^{-1})$ st"uckweise stetig,
besitze h"ochstens abz"ahlbar viele Sprungstellen $s_k \in (0,\infty)$,
die sich nirgends im Endlichen h"aufen,
und $p(\cdot)$ sei zwischen diesen Spr"ungen stetig differenzierbar. \\
Ferne gelte:
\begin{enumerate}
\item[(a)] Das Tripel $\big(x_*(\cdot),u_*(\cdot),p(\cdot)\big)$
           erf"ullt (\ref{SatzWMPAufgabeZB1})--(\ref{SatzWMPAufgabeZB3}) mit $\lambda_0=1$ in Theorem \ref{SatzWMPAufgabeZB}.        
\item[(b)] F"ur jedes $t \in \R_+$ ist die Funktion $H\big(t,x,u,p(t),1\big)$ konkav in $(x,u)$
           und es sind die Funktionen $g_j(t,x)$, $j=1,...,l$, konvex bez"uglich $x$ auf $U_\gamma$.
\end{enumerate}
Dann ist $\big(x_*(\cdot),u_*(\cdot)\big)$ ein schwaches lokales Minimum der Aufgabe (\ref{WMPAufgabe1})--(\ref{WMPAufgabe4}).
\end{theorem}

\begin{beispiel}[Abbau einer nicht erneuerbaren Ressource] \label{ExampleRessource}
{\textrm Wir betrachten die Aufgabe
\begin{eqnarray}
&& \label{Ressource1} J\big(x(\cdot),y(\cdot),u(\cdot)\big)
   =\int_0^\infty e^{-\varrho t}\big[pf\big(u(t)\big)-ry(t)-qu(t)\big] \, dt \to \sup, \\
&& \label{Ressource2} \dot{x}(t) = -u(t),\quad \dot{y}(t)=cf\big(u(t)\big), \quad x(0)=x_0>0,\quad y(0)=y_0\geq 0, \\
&& \label{Ressource3} x(t) \geq 0, \qquad u(t) \geq 0, \qquad b,c,\varrho, q,r >0, \qquad \varrho - rc>0.
\end{eqnarray}
Die Funktion $f$ sei zweimal stetig differenzierbar, $f'>0$, $f'(0)<\infty$, $f''<0$ und es sei
$f'(u) \to 0$ f"ur $u \to \infty$.
In der vorliegenden Formulierung der Aufgabe wurde im Vergleich zu Seierstad \& Syds\ae ter \cite{Seierstad} die
Restriktion $\liminf\limits_{t \to \infty} x(t) \geq 0$ durch die Zustandsbeschr"ankung $x(t) \geq 0$ in (\ref{Ressource3})
ersetzt. \\[2mm]
"Okonomische Interpretation:
$x(t)$ bezeichnet die Menge einer nat"urlichen Ressource und $u(t)$ ist die industrielle Abbaurate dieser Ressource.
Auf Basis der Ressource werden G"uter mit der Produktionsrate $f\big(u(t)\big)$ hergestellt.
Die Kosten der Herstellung einer Produktionseinheit ist $q$ und der Preis einer G"utereinheit am Markt betr"agt $p$.
Bei der Herstellung der G"uter entstehen proportional zur Produktion Abf"alle,
deren Gesamtmenge durch $y(t)$ beschrieben wird.
Die Kosten der Beseitigung der negativen Auswirkungen der Abfallmenge sind $ry(t)$.
Im Weiteren gehen wir von dem Preis $p=1$ aus. \\[2mm]
Wir pr"ufen die Voraussetzungen an die Aufgabe:
Mit der Festlegung $L(t)=C_0e^{at}$ und $\nu(t)=e^{-at}$ mit $0<a<\varrho$ gen"ugt die Aufgabe den Voraussetzungen (A$_1$)--(A$_3$).
Weiterhin sind die Dynamiken
$$\varphi_1(t,x,y,u)=-u, \qquad \varphi_2(t,x,y,u)=cf(u)$$
unabh"angig von den Zustandsvariablen $x,y$ und die Bedingungen (\ref{WMPBedingungZB}), (\ref{WMPBedingungZB2}) sind erf"ullt.
Wir stellen die Optimalit"atsbedingungen von Theorem \ref{SatzWMPAufgabeZB} mit $\lambda_0=1$ auf:
\begin{enumerate}
\item[(a)] Die Pontrjagin-Funktion der Aufgabe (\ref{Ressource1})--(\ref{Ressource3}) lautet
           $$H(t,x,y,u,p_1,p_2,1) = p_1 (-u)+p_2 cf(u) + e^{-\varrho t}[f(u)-ry-qu].$$
\item[(b)] Die Adjungierten gen"ugen den Gleichungen
           $$p_1(t)=\int_t^\infty e^{-as} \, d\mu(s), \qquad \dot{p}_2(t)=r e^{-\varrho t} \Rightarrow p_2(t)=-\frac{r}{\varrho}e^{-\varrho t} + K.$$
           Das auf der Menge $T=\{t \in \overline{\R}_+ \,|\, x_*(t)=0\}$ konzentrierte regul"are Ma"s $\mu$ ist nichtnegativ.
           Daher ist $p_1(t) \geq 0$ "uber $\R_+$ und monoton fallend.
           Ferner erhalten wir $K=0$ aus den Transversalit"atsbedingungen bez"uglich dem Zustand $y$.
\item[(c)] Die Maximumbedingung k"onnen wir auf folgende Aufgabe reduzieren
           $$\max_{u \geq 0} \Big[ -p_1(t) u +c p_2(t) f(u) + e^{-\varrho t}[f(u)-qu]\Big].$$
           Das Einsetzen der Darstellung f"ur $p_2(t)$ liefert weiterhin mit $d=(\varrho -rc)/\varrho$:
           $$\max_{u \geq 0} \Big( d f(u)e^{-\varrho t}-u\big(p_1(t)+ qe^{-\varrho t}\big)\Big).$$
\end{enumerate}
Die Reduktion der Maximumbedingung f"uhrt f"ur festes $t$ zu der Funktion
$$g(u)= d f(u)e^{-\varrho t}-u\big(p_1(t)+ qe^{-\varrho t}\big).$$
Diese Funktion ist zweimal stetig differenzierbar und es gilt
$$g'(u)= \big(d f'(u)-q\big)e^{-\varrho t}-p_1(t), \quad g''(u)=df''(u)e^{-\varrho t}, \quad d=\frac{\varrho -rc}{\varrho}>0.$$
Daher ist $g$ streng konkav und besitzt auf der Menge $U=\{u\geq 0\}$ ein Maximum,
da $f'(u)>0$ und $f'(u) \to 0$ f"ur $u \to \infty$ gelten.
Wir diskutieren drei F"alle:
\begin{enumerate}
\item[(A)] $df'(0)\leq q$: In diesem Fall ist $g'(0) \leq 0$ und man erh"alt
           $$u_*(t) \equiv 0, \quad x_*(t) \equiv x_0, \quad
             y_*(t)=y_0 + cf(0)t, \quad p_1(t) \equiv 0, \quad p_2(t)=-\frac{r}{\varrho}e^{-\varrho t}.$$
           Da die Zustandsbeschr"ankung nichtaktiv ist, gelten die Voraussetzungen und Optimalit"atsbedingungen aus Theorem \ref{SatzWMPAufgabeZB} 
           f"ur $\nu(t)=e^{-at}$ mit $0<a<\varrho$.
\item[(B)] $df'(0)> q$ und $p_1(0)=0$:
           Aus $g'(u)=\big(d f'(u)-q\big)e^{-\varrho t}=0$ erhalten wir die optimale Strategie $u_*(t)=u_0>0$ f"ur alle $t \in \R_+$.
           Also gilt $x_*(t)=x_0-u_0t$ auf $\R_+$, was der Zustandsbeschr"ankung widerspricht.
\item[(C)] $df'(0)> q$ und $p_1(0)>0$:
           Wegen $p_1(0)>0$ wird die Ressource vollst"andig abgebaut.
           Andernfalls w"are $p_1(t)=p_1(0)>0$ "uber $\R_+$, was (\ref{SatzWMPAufgabeZB2}) widerspricht.
           Da die Ressource vollst"andig abgebaut wird,
           gibt es ein $t'>0$ mit $x_*(t)>0$ f"ur $t \in [0,t')$ und $x_*(t)=0$ f"ur $t\geq t'$.
           Demnach folgt unmittelbar $u_*(t)=0$ f"ur $t\geq t'$. \\
           F"ur $t\geq t'$ ist $p_1(\cdot)$ monoton fallend.
           Ferner erhalten wir f"ur $t \in \R_+$ die Beziehung
           $$g'(u)=0 \qquad\Rightarrow\qquad f'\big(u(t)\big)=\frac{1}{d}(q+p_1(t)e^{\varrho t}).$$
           W"urde demnach die Adjungierte $p_1(\cdot)$ f"ur $t \geq t'$ eine Unstetigkeitstelle besitzen,
           dann folgt aus der Monotonie von $p_1(\cdot)$, dass die Abbaurate sich wieder sprunghaft vergr"o"sert.
           Diese Steuerung f"uhrt zu einem erneuten Abbau der Ressource, obwohl diese bereits vollst"andig aufgebraucht ist.
           Daher ist die Adjungierte stetig. \\
           F"ur die Adjungierte erhalten wir damit
           $$p_1(t) = \big(df'(0)-q\big)e^{-\varrho t'} \mbox{ f"ur } t \leq t', \qquad
             p_1(t) = \big(df'(0)-q\big)e^{-\varrho t} \mbox{ f"ur } t \geq t'.$$
           Wir zeigen noch, dass der Zeitpunkt $t'$ existiert und eindeutig ist:
           Durch
           $$f'\big(u_\tau(t)\big)=\frac{1}{d}(q+p_1(0)e^{\varrho t})=\frac{1}{d}(q+[df'(0)-q]e^{\varrho(t-\tau)}), \quad t \in [0,\tau],$$
           und $u_\tau(t)=0$ f"ur $t \geq \tau$ wird wegen $f'\big(u_\tau(\tau)\big)=f'(0)$ eine Familie $u_\tau(\cdot)$ stetiger Funktionen definiert.
           Dabei gilt $f'\big(u_\tau(t)\big) < f'\big(u_s(t)\big)$, d.\,h. $u_\tau(t) > u_s(t)$, f"ur alle $t \in [0,\tau]$ und $\tau>s$.
           Damit ist die Familie
           $$U(\tau):= \int_0^\infty u_\tau(t) \, dt$$
           streng monoton wachsend und es gelten $U(0)=0$, $U(\tau) \to \infty$ f"ur $\tau \to \infty$.
           Der Parameter $t'$ ergibt sich dann aus der Bedingung $U(t')=x_0$. \hfill $\square$
\end{enumerate}
Die Funktion $H\big(t,x,u,p_1(t),p_2(t),1\big) = p_1(t) (-u)+p_2(t) cf(u) + e^{-\varrho t}[f(u)-ry-qu]$ ist nach den Voraussetzungen an die Funktion $f(u)$ konkav.
Damit liefern die ermittelten Kandidaten in den F"allen (A) und (C) schwache lokale Minimalstellen der Aufgabe. \hfill $\square$}
\end{beispiel}

\lhead[\thepage \hspace*{1mm} Starkes lokales Minimum]{}
\rhead[]{Starkes lokales Minimum \hspace*{1mm} \thepage}
\section{Starkes lokales Minimum "uber unendlichem Zeithorizont} \label{KapitelStrong}
\subsection{Die Aufgabenstellung}
In diesem Kapitel betrachten wir starke lokale Minimalstellen der Aufgabe
\begin{eqnarray}
&& \label{PAUH1} J\big(x(\cdot),u(\cdot)\big) = \int_0^\infty \omega(t)f\big(t,x(t),u(t)\big) \, dt \to \inf, \\
&& \label{PAUH2} \dot{x}(t) = \varphi\big(t,x(t),u(t)\big), \\
&& \label{PAUH3} h_0\big(x(0)\big)=0, \qquad \lim_{t \to \infty} h_1\big(t,x(t)\big)=0, \\
&& \label{PAUH4} u(t) \in U \subseteq \R^m, \quad U \not= \emptyset, \\
&& \label{PAUH5} g_j\big(t,x(t)\big) \leq 0 \quad \mbox{f"ur alle } t \in \R_+, \quad j=1,...,l.
\end{eqnarray}
Dabei ist $\omega(\cdot) \in L_1(\R_+,\R_+)$ und es gelten f"ur die eingehenden Abbildungen
$$f:\R \times \R^n \times \R^m \to \R, \qquad \varphi:\R \times \R^n \times \R^m \to \R^n, \qquad g_j: \R \times \R^n \to \R,$$
sowie f"ur die Randbedingungen
$$h_0:\R^n \to \R^{s_0}, \qquad h_1:\R \times \R^n \to \R^{s_1}.$$
Wir nennen die Trajektorie $x(\cdot)$ eine L"osung des dynamischen Systems (\ref{PAUH2}),
falls $x(\cdot)$ auf $\R_+$ definiert ist und auf jedem endlichen Intervall die Dynamik mit Steuerung $u(\cdot)$
im Sinn von Carath\'eodory l"ost. \\[2mm]
Zu $x(\cdot)$ bezeichne $V_\gamma$ die Menge
$V_\gamma= \{ (t,x) \in \overline{\R}_+ \times \R^n \,|\, \|x-x(t)\| \leq \gamma\}$.
Dann geh"oren zur Menge $\mathscr{B}_{\textrm Lip}$ diejenigen
$x(\cdot) \in W^1_\infty(\R_+,\R^n)$,
f"ur die es zu jeder kompakten Menge $U_1 \subseteq \R^m$ eine Zahl $\gamma>0$ derart gibt,
dass auf der Menge $V_\gamma \times U_1$ die Abbildungen
\begin{enumerate}
\item[(B)] $f(t,x,u)$, $\varphi(t,x,u)$, $g_j(t,x)$ und $h_0(x)$, $h_1(t,x)$ gleichm"a"sig stetig und
           gleichm"a"sig stetig differenzierbar bez"uglich $x$ sind.
\end{enumerate}

Der Steuerungsprozess $\big(x(\cdot),u(\cdot)\big) \in W^1_\infty(\R_+,\R^n;\nu) \times L_\infty(\R_+,U)$
hei"st zul"assig in der Aufgabe (\ref{PAUH1})--(\ref{PAUH5}),
falls $\big(x(\cdot),u(\cdot)\big)$ dem System (\ref{PAUH2}) gen"ugt,
die Randbedingungen (\ref{PAUH3}) und Restriktionen (\ref{PAUH4}) erf"ullt,
sowie das Lebesgue-Integral in (\ref{PAUH1}) endlich ist.
Die Menge $\mathscr{B}_{\textrm adm}$ bezeichnet die Menge der zul"assigen Steuerungsprozesse der Aufgabe (\ref{PAUH1})--(\ref{PAUH5}). \\
Der zul"assige Steuerungprozess $\big(x_*(\cdot),u_*(\cdot)\big)$ ist eine starke lokale Minimalstelle\index{Minimum, schwaches lokales!--, starkes lokales} in
der Aufgabe (\ref{PAUH1})--(\ref{PAUH5}),
falls eine Zahl $\varepsilon > 0$ derart existiert, dass die Ungleichung
$$J\big(x(\cdot),u(\cdot)\big) \geq J\big(x_*(\cdot),u_*(\cdot)\big)$$
f"ur alle $\big(x(\cdot),u(\cdot)\big) \in \mathscr{B}_{\textrm adm}$ mit 
$\|x(\cdot)-x_*(\cdot)\|_\infty \leq \varepsilon$ gilt.
       
       \subsection{Das Pontrjaginsche Maximumprinzip f"ur eine Grundaufgabe} \label{AbschnittPMPUH}
Im Weiteren bezeichnet $H: \R \times \R^n \times \R^m \times \R^n \times \R \to \R$ die Pontrjagin-Funktion
$$H(t,x,u,p,\lambda_0) = \langle p, \varphi(t,x,u) \rangle-\lambda_0 \omega(t)f(t,x,u).$$

\begin{theorem}[Pontrjaginsches Maximumprinzip] \label{SatzPAUH}\index{Pontrjaginsches Maximumprinzip}
Es sei $\big(x_*(\cdot),u_*(\cdot)\big) \in \mathscr{B}_{\textrm adm} \cap \mathscr{B}_{\textrm Lip}$.
Weiterhin nehmen wir an,
dass
\begin{equation} \label{PMPBedingung}
\int_0^\infty \big\|\varphi\big(t,x_*(t),u_*(t)\big)\big\| \, dt < \infty, \qquad \int_0^\infty \big\|\varphi_x\big(t,x_*(t),u_*(t)\big)\big\| \, dt < \infty
\end{equation}
ausfallen und es m"oge zu jedem $\delta>0$ ein $T>0$ existieren mit
\begin{eqnarray}
&& \int_T^\infty \big\| \varphi\big(t,x(t),u_*(t)\big)-\varphi\big(t,x'(t),u_*(t)\big) - \varphi_x\big(t,x_*(t),u_*(t)\big)\big(x(t)-x'(t)\big) \big\| \, dt
   \nonumber \\
&& \label{PMPBedingung2} \hspace*{20mm} \leq \delta \|x(\cdot)-x'(\cdot)\|_\infty
\end{eqnarray}
f"ur alle $x(\cdot), x'(\cdot) \in W^1_\infty(\R_+,\R^n)$ mit $\|x(\cdot)-x_*(\cdot)\|_\infty < \gamma$, $\|x'(\cdot)-x_*(\cdot)\|_\infty < \gamma$. \\[2mm]
Ist $\big(x_*(\cdot),u_*(\cdot)\big)$ ein starkes lokales Minimum der Aufgabe (\ref{PAUH1})--(\ref{PAUH4}),
dann existieren nicht gleichzeitig verschwindende Multiplikatoren $\lambda_0 \geq 0$,
$p(\cdot) \in W^1_\infty(\R_+,\R^n)$ und $l_i \in \R^{s_i}$, $i=0,1$, derart, dass
\begin{enumerate}
\item[(a)] die Funktion $p(\cdot)$ fast "uberall der adjungierten Gleichung\index{adjungierte Gleichung}
           \begin{equation}\label{SatzPAUH1}
           \dot{p}(t) = -\varphi_x^T\big(t,x_*(t),u_*(t)\big) p(t) + \lambda_0 \omega(t)f_x\big(t,x_*(t),u_*(t)\big)
           \end{equation}
           gen"ugt und die Transversalit"atsbedingungen\index{Transversalit"atsbedingungen}
           \begin{equation}\label{SatzPAUH2}
           p(0) = {h_0'}^T\big(x_*(0)\big)l_0, \qquad \lim_{t \to \infty} p(t)= - \lim_{t \to \infty}h_{1x}^T\big(t,x_*(t)\big) l_1
           \end{equation}
           erf"ullt;
\item[(b)] in fast allen Punkten $t \in \R_+$ die Maximumbedingung gilt:
           \begin{equation}\label{SatzPAUH3}
           H\big(t,x_*(t),u_*(t),p(t),\lambda_0\big) = \max_{u \in U} H\big(t,x_*(t),u,p(t),\lambda_0\big).
           \end{equation}
\end{enumerate}
\end{theorem}

\begin{beispiel} {\textrm Wir betrachten die Aufgabe
\begin{eqnarray*}
&& J\big(x(\cdot),u(\cdot)\big) = \int_0^\infty e^{-\varrho t} \big(1-u(t)\big) x(t) \, dt \to \sup,\\
&& \dot{x}(t)=u(t)x(t), \quad x(0)=1, \quad \lim_{t \to \infty} x(t)=x_1>1, \quad u(t) \in [0,1], \quad \varrho \in (0,1).
\end{eqnarray*}
Wir stellen zun"achst die Pontrjagin-Funktion auf:
$$H(t,x,z,u,p,q,1) = pux+\lambda_0 e^{-\varrho t}(1-u)x.$$
Mit Hilfe der Bedingungen (\ref{SatzPAUH1})--(\ref{SatzPAUH3}) k"onnen wir den Fall $\lambda_0=0$ ausschlie"sen.
Weiterhin ergeben sich der Steuerungsprozess
$$x_*(t) = \left\{ \begin{array}{ll} e^t, & t \in [0,\tau), \\ x_1, & t \in [\tau, \infty), \end{array}\right. \quad
  u_*(t)= \left\{ \begin{array}{ll} 1, & t \in [0,\tau), \\ 0, & t \in [\tau, \infty), \end{array}\right. \quad \tau=\ln x_1$$
und die Adjungierte
$$p(t) = \left\{ \begin{array}{ll}
          e^{(1-\varrho)\tau} e^{-t}, & t \in [0,\tau), \\
          \frac{\varrho-1}{\varrho}e^{-\varrho \tau} + \frac{1}{\varrho} e^{-\varrho t}, & t \in [\tau, \infty). \end{array}\right.$$
F"ur die Adjungierte gilt dabei im Unendlichen $\displaystyle\lim_{t \to \infty} p(t)= \frac{\varrho-1}{\varrho}e^{-\varrho \tau} \not=0$.
Damit konnten wir aus den notwendigen Bedingungen (\ref{SatzPAUH1})--(\ref{SatzPAUH3}) des Maximumprinzips einen eindeutigen
Kandidaten bestimmen. \hfill $\square$}
\end{beispiel}

\begin{beispiel} \label{BeipielRWUnendlich}
{\textrm Wir betrachten die Aufgabe
\begin{eqnarray*}
&& J\big(x(\cdot),u(\cdot)\big)=\int_0^\infty e^{-\varrho t}\big(1-u(t)\big)x(t) \, dt \to \sup, \\
&& \dot{x}(t)=u(t)x(t), \qquad x(0)=1, \qquad u \in [0,1], \qquad \varrho \in (0,1)
\end{eqnarray*}
mit der Budgetbeschr"ankung
$$\int_0^\infty e^{-\varrho t}x(t) \, dt = Z, \qquad Z> \frac{1}{\varrho}.$$
Bez"uglich der Budgetbeschr"ankung f"uhren wir die folgende Zustandsgleichung mit Randwert im Unendlichen ein:
$$\dot{z}(t)=e^{-\varrho t}x(t), \qquad z(0)=0, \qquad \lim_{t \to \infty} z(t) = Z.$$
Die Zustandsgleichung und -beschr"ankung f"ur die Trajektorie $z(\cdot)$ ergibt sich aus der isoperimetrischen Nebenbedingung in Form einer Budgetbeschr"ankung
$$\int_0^\infty e^{-\varrho t}x(t) \, dt \leq Z.$$
Offensichtlich ist $\dot{z}(t) >0$ auf $\R_+$ und damit $z(t)$ streng monoton wachsend.
Demzufolge kann die Beschr"ankung $z(t) \leq Z$ erst im Unendlichen aktiv werden und greift nur durch das Verhalten in $t=\infty$ in die
gestellte Aufgabe ein.
Da stets $\dot{x}(t)\geq 0$ ist,
muss f"ur jede zul"assige Trajektorie $e^{-\varrho t}x(t) \to 0$ f"ur $t \to \infty$ gelten,
denn nur dann ist $z(t) \leq Z$ erf"ullt.
Damit erhalten wir f"ur zul"assige Steuerungsprozesse zun"achst im Zielfunktional
\begin{eqnarray*}
    J\big(x(\cdot),z(\cdot),u(\cdot)\big)
&=& \int_0^\infty e^{-\varrho t}\big(1-u(t)\big)x(t) \, dt = \int_0^\infty \dot{z}(t) \, dt - \int_0^\infty e^{-\varrho t} \dot{x}(t) \, dt \\
&=& \int_0^\infty \dot{z}(t) \, dt + 1 - \varrho \int_0^\infty e^{-\varrho t} x(t) \, dt \leq 1 + (1-\varrho) Z.
\end{eqnarray*}
Es ergibt sich also f"ur das Zielfunktional die obere Schranke $1 + (1-\varrho) Z$.
D.\,h., dass jeder Steuerungsprozess $\big(x(\cdot),z(\cdot),u(\cdot)\big)$,
f"ur den die Zustandsbeschr"ankung $z(t) \leq Z$ im Unendlichen aktiv wird,
global optimal ist und $J\big(x(\cdot),z(\cdot),u(\cdot)\big) = 1 + (1-\varrho) Z$ gilt. \\
Die Voraussetzungen des Pontrjaginschen Maximumprinzips an einen zul"assigen Steuerungsprozess sind in dem vorliegenden Beispiel genau dann erf"ullt,
wenn die Steuerung $u(\cdot)$ dem Raum $L_1(\R_+,[0,1])$ angeh"ort.
Denn in diesem Fall gelten $x(\cdot) \in C_{\lim}(\R_+,\R)$ f"ur die korrespondierende Trajektorie und die Bedingungen 
(\ref{PMPBedingung}), (\ref{PMPBedingung2}).
Wir diskutieren zwei optimale Steuerungsprozesse:
\begin{enumerate}
\item[(A)] In diesem Beispiel liefert der Steuerungsprozess
           $$x_*(t) = e^{\alpha t},\qquad z_*(t)= \frac{1}{\alpha - \varrho} (e^{(\alpha-\varrho)t}-1),\qquad
             u_*(t)=\alpha,\qquad \alpha=\varrho-\frac{1}{Z} \in (0, \varrho)$$
           ein globales Maximum.
           Da die vorgeschlagene Steuerung $u_*(\cdot)$ "uber $\R_+$ nicht integrierbar ist,
           gelten weder $x_*(\cdot) \in C_{\lim}(\R_+,\R)$ noch die Bedingungen (\ref{PMPBedingung}). \\
           Das Maximumprinzip ist auf diesen Steuerungsprozess nicht anwendbar.
\item[(B)] Ebenfalls stellt der Steuerungsprozess
           $$y_*(t) = \left\{ \begin{array}{ll} e^t,& t \in [0,\tau), \\ e^\tau, & t \in [\tau,\infty), \end{array} \right. \quad
             w_*(t) = \left\{ \begin{array}{ll} 1,& t \in [0,\tau), \\ 0, & t \in [\tau,\infty) \end{array} \right.$$
           mit dem Umschaltzeitpunkt $\tau >0$, der der Bedingung
           $$e^{(1-\varrho)\tau}\bigg(\frac{1}{\varrho}+\frac{1}{1-\varrho}\bigg) = Z+\frac{1}{1-\varrho}, \qquad Z> \frac{1}{\varrho},$$
           gen"ugt,
           ein globales Maximum dar.
           Die zugeh"orige Trajektorie $z_*(\cdot)$ lautet
           $$z_*(t) = \left\{ \begin{array}{ll} \frac{1}{1-\varrho}\big(e^{(1-\varrho)t} - 1 \big) ,& t \in [0,\tau), \\
                                                z(\tau) + \frac{1}{\varrho}\big(e^{(1-\varrho)\tau} - e^{\tau-\varrho t} \big), & t \in [\tau,\infty).
             \end{array} \right.$$
           Da die Steuerung $w_*(\cdot)$ dem Raum $L_1(\R_+,[0,1])$ angeh"ort,
           gelten s"amtliche Voraussetzungen von Theorem \ref{SatzPAUH}.
           Mit den Multiplikatoren
           $$\lambda_0=1, \qquad p(t)= e^{-\varrho t}, \qquad q(t)= \varrho-1$$
           sind dann die notwendigen Bedingungen (\ref{SatzPAUH1})--(\ref{SatzPAUH3}) erf"ullt. \hfill $\square$
\end{enumerate}}
\end{beispiel}

F"ur die Aufgabe (\ref{PAUH1})--(\ref{PAUH4}) lassen sich Aussagen "uber die Normalform des Pontrjaginschen Maximumprinzips und zu diversen
Transversalit"atsbedingungen ableiten. \\[2mm]
Wir betrachten zun"achst die Aufgabe (\ref{PAUH1})--(\ref{PAUH4}) mit freiem rechten Endpunkt im Unendlichen.
Dann gen"ugt die Adjungierte $p(\cdot)$ nach Theorem \ref{SatzPAUH} dem Randwertproblem
$$\dot{p}(t) = -\varphi_x^T\big(t,x_*(t),u_*(t)\big) p(t) + \lambda_0 \omega(t)f_x\big(t,x_*(t),u_*(t)\big), \qquad
  \lim_{t \to \infty} p(t)=0.$$
  
\begin{folgerung} \label{FolgerungPAUH1}
In der Aufgabe (\ref{PAUH1})--(\ref{PAUH4}) mit freiem rechten Endpunkt im Unendlichen
ergeben sich aus dem Maximumprinzip unmittelbar die ``nat"urlichen'' Transversalit"atsbedingungen:
$$\lim_{t \to \infty} p(t) =0, \qquad \lim_{t \to \infty} \langle p(t),x(t) \rangle = 0 \mbox{ f"ur alle } x(\cdot) \in W^1_\infty(\R_+,\R^n).$$
\end{folgerung}

Nach Voraussetzung des Maximumprinzips ist
$$\int_0^\infty \big\|\varphi_x\big(t,x_*(t),u_*(t)\big)\big\| \, dt < \infty.$$
Unter der Annahme $\lambda_0=0$ w"urde die Adjungierte als L"osung der Gleichung
$$\dot{p}(t) = -\varphi_x^T\big(t,x_*(t),u_*(t)\big) p(t), \qquad \lim_{t \to \infty} p(t)=0,$$
im Widerspruch zu Theorem \ref{SatzPAUH} identisch verschwinden.

\begin{folgerung} \label{FolgerungPAUH2}
In der Aufgabe (\ref{PAUH1})--(\ref{PAUH4}) mit freiem rechten Endpunkt im Unendlichen gilt $\lambda_0 \not= 0$ und wir k"onnen ohne Einschr"ankung
$\lambda_0=1$ annehmen.
\end{folgerung}

Wegen der Integrierbarkeit der Abbildung $t \to \varphi_x\big(t,x_*(t),u_*(t)\big)$ "uber $\R_+$ k"onnen wir die 
die in $t=0$ normalisierten Fundamentalmatrizen $Y_*(t)$ bzw. $Z_*(t)$ der homogenen Systeme
$$\dot{y}(t)=\varphi_x\big(t,x_*(t),u_*(t)\big) y(t), \qquad \dot{z}(t)=-\varphi_x^T\big(t,x_*(t),u_*(t)\big) z(t)$$
im Rahmen des Raumes $C_{\lim}(\R_+,\R^n)$ betrachten.
Es bezeichne ferner $y_\xi(\cdot) \in C_{\lim}(\R_+,\R^n)$ die Funktion $y_\xi(t)=Y_*(t) Y^{-1}_*(T) \xi$ mit $\xi \in \R^n$ und $\|\xi\|=1$.
Dann ergibt sich auf die gleiche Weise wie in Aseev \& Kryazhimskii und Aseev \& Veliov \cite{AseKry,AseVel,AseVel2,AseVel3} oder
Tauchnitz \cite{TauchnitzWMPIHOC,TauchnitzPMPIHOC,TauchnitzOC} die Beziehung
$$\langle p(t), y_\xi(t) \rangle = \Big\langle  p(T) + Z_*(T)\int_T^t \omega(s) Z_*^{-1}(s)f_x\big(s,x_*(s),u_*(s)\big) ds , \xi \Big\rangle.$$
Unter Verwendung der ``nat"urlichen'' Transversalit"atsbedingungen in Folgerung \ref{FolgerungPAUH1} und wegen der Willk"urlichkeit von $\xi$
erhalten wir daraus die Integraldarstellung
der Arbeiten von Aseev \& Kryazhimskii und Aseev \& Veliov \cite{AseKry,AseVel,AseVel2,AseVel3}:

\begin{folgerung} \label{FolgerungPAUH3}
Es gen"uge $\big(x_*(\cdot),u_*(\cdot)\big)$ den Voraussetzungen des Pontrjaginschen Maximumprinzips \ref{SatzPAUH}.
Ist $\big(x_*(\cdot),u_*(\cdot)\big)$ ein starkes lokales Minimum der Aufgabe (\ref{PAUH1})--(\ref{PAUH4}) mit freiem rechten Endpunkt im Unendlichen,
dann besitzt die Adjungierte $p(\cdot)$ die Darstellung\index{Adjungierte!eindeutig@--, eindeutige Darstellung}
$$p(t)= -Z_*(t) \int_t^\infty \omega(s) Z^{-1}_*(s) f_x\big(s,x_*(s),u_*(s)\big) \, ds.$$
Dabei ist $Z_*(t)$ die in $t=0$ normalisierte Fundamentalmatrix des linearen Systems
$$\dot{z}(t)=-\varphi^T_x\big(t,x_*(t),u_*(t)\big) z(t).$$
\end{folgerung}

In der Aufgabe (\ref{PAUH1})--(\ref{PAUH4}) seien nun gewisse Randwerte im Unendlichen explizit gegeben, d.\,h. $h_1\big(t,x(t)\big)=x(t)-x_1$.
Wir schlie"sen dabei nicht aus,
dass dabei gewisse Komponenten von $x_1$ nicht fest vorgegeben, sondern ohne Einschr"ankung sind.
Wir sprechen dabei von expliziten Randbedingungen wenn f"ur $x(t)=\big(x_1(t),...,x_n(t)\big)$ gilt:
$$\lim_{t \to \infty} x_i(t)=x_i \in \R, \qquad \lim_{t \to \infty} x_j(t) \mbox{ frei}, \qquad i=1,...,l,\; j=l+1,...,n.$$
Daher verschwinden einerseits bei expliziten Randwerten im Unendlichen die entsprechenden Komponenten $x_i(t)-x_{i*}(t)$ f"ur $t \to \infty$, $i=1,...,l$.
F"ur diejenigen Komponenten, f"ur die die Randwerte im Unendlichen frei sind,
verschwinden die entsprechenden Komponenten $p_j(t)$ der Adjungierten f"ur $t \to \infty$, $j=l+1,...,n$.
Damit k"onnen wir festhalten:

\begin{folgerung} \label{FolgerungPAUH4}
In der Aufgabe (\ref{PAUH1})--(\ref{PAUH4}) mit expliziten Randbedingungen im Unendlichen
ergibt sich aus dem Maximumprinzip unmittelbar die ``nat"urliche'' Transversalit"atsbedingung\index{Transversalit"atsbedingungen!nat@--, nat"urliche}:
$$\lim_{t \to \infty} \langle p(t),x(t)-x_*(t) \rangle = 0 \quad \mbox{ f"ur alle zul"assigen } x(\cdot) \in W^1_\infty(\R_+,\R^n).$$
\end{folgerung}

In der Aufgabe (\ref{PAUH1})--(\ref{PAUH4}) mit freiem rechten Endpunkt verschwinde die Verteilungsfunktion $\omega(\cdot) \in L_1(\R_+,\R_+)$
im Unendlichen.
Dann folgen aus der gleichm"a"sigen Stetigkeit der Abbildungen $f,\,\varphi$ und aus der Beschr"anktheit des Steuerungsprozesses $\big(x_*(\cdot),u_*(\cdot)\big)$ 
die Grenzwerte
$$\lim_{t \to \infty} \omega(t) f\big(t,x_*(t),u_*(t)\big)=0, \qquad
  \lim_{t \to \infty} \big\langle p(t)\,,\, \varphi\big(t,x_*(t),u_*(t)\big)\big\rangle=0.$$

\begin{folgerung} \label{FolgerungPAUH5}
In der Aufgabe (\ref{PAUH1})--(\ref{PAUH4}) mit freiem Endpunkt im Unendlichen besitze die Verteilungsfunktion $\omega(\cdot)$ einen Grenzwert im
Unendlichen, d.\,h.
$$\lim_{t \to \infty} \omega(t)=0.$$
Dann ergibt sich die Bedingung von Michel \index{Transversalit"atsbedingungen!von@-- von Michel}:
$$\lim_{t \to \infty} H\big(t,x_*(t),u_*(t),p(t),\lambda_0\big)=0.$$
\end{folgerung}

Es bezeichnet $\mathscr{H}(t,x,p) = \sup\limits_{u \in U} H(t,x,u,p,1)$ die Hamilton-Funktion $\mathscr{H}$ im normalen Fall.

\begin{theorem} \label{SatzHBPMPUH}
In der Aufgabe (\ref{PAUH1})--(\ref{PAUH4}) mit expliziten Randbedingungen
sei $\big(x_*(\cdot),u_*(\cdot)\big) \in \mathscr{B}_{\textrm Lip} \cap \mathscr{B}_{\textrm adm}$ 
und es sei $p(\cdot):\R_+ \to \R^n$. Ferner gelte:
\begin{enumerate}
\item[(a)] Das Tripel $\big(x_*(\cdot),u_*(\cdot),p(\cdot)\big)$
           erf"ullt (\ref{SatzPAUH1})--(\ref{SatzPAUH3}) in Theorem \ref{SatzPAUH} mit $\lambda_0=1$.        
\item[(b)] F"ur jedes $t \in \R_+$ ist die Funktion $\mathscr{H}\big(t,x,p(t)\big)$ konkav in $x$ auf $V_\gamma(t)$.
\end{enumerate}
Dann ist $\big(x_*(\cdot),u_*(\cdot)\big)$ ein starkes lokales Minimum der Aufgabe (\ref{PAUH1})--(\ref{PAUH4}) mit expliziten Randbedingungen.
\end{theorem}

\begin{beispiel} {\textrm Wir betrachten die Aufgabe
\begin{eqnarray*}
&& J\big(x(\cdot),u(\cdot)\big) = \int_0^\infty e^{-\varrho t} \big(1-u(t)\big) x(t) \, dt \to \sup,\\
&& \dot{x}(t)=u(t)x(t), \quad x(0)=1, \quad \lim_{t \to \infty} x(t)=x_1>1, \quad u(t) \in [0,1], \quad \varrho \in (0,1).
\end{eqnarray*}
Aus den notwendigen Bedingungen (\ref{SatzPAUH1})--(\ref{SatzPAUH3}) ergeben sich
\begin{eqnarray*}
x_*(t) &=& \left\{ \begin{array}{ll} e^t, & t \in [0,\tau), \\ x_1, & t \in [\tau, \infty), \end{array}\right. \quad
  u_*(t)= \left\{ \begin{array}{ll} 1, & t \in [0,\tau), \\ 0, & t \in [\tau, \infty), \end{array}\right. \quad \tau=\ln x_1, \\
p(t) &=& \left\{ \begin{array}{ll}
        e^{(1-\varrho)\tau} e^{-t}, & t \in [0,\tau), \\
        \frac{\varrho-1}{\varrho}e^{-\varrho \tau} + \frac{1}{\varrho} e^{-\varrho t}, & t \in [\tau, \infty).
        \end{array}\right.
\end{eqnarray*}
Offenbar ist die Hamilton-Funktion $\mathscr{H}$ konkav in $x$ und damit $\big(x_*(\cdot),u_*(\cdot)\big)$
ein starkes lokales Minimum der Aufgabe.  \hfill $\square$}
\end{beispiel}

\begin{beispiel} {\textrm F"ur den linear-quadratische Regler
\begin{eqnarray*}
&& J\big(x(\cdot),u(\cdot)\big) = \int_0^\infty e^{-2t} \cdot \frac{1}{2}\big( x^2(t)+u^2(t)\big) \, dt \to \inf, \\
&& \dot{x}(t) = 2 x(t)+u(t), \qquad x(0)=2, \qquad u(t) \in \R
\end{eqnarray*}
liefern die Bedingungen (\ref{SatzPAUH1})--(\ref{SatzPAUH3}) den Steuerungsprozess und die Adjungierte
$$x_*(t)=2e^{(1-\sqrt{2})t}, \quad u_*(t)=-2(1+\sqrt{2})e^{(1-\sqrt{2})t}, \quad p(t)=e^{-2t}u_*(t).$$
Die Hamilton-Funktion $\mathscr{H}$ ist konkav bez"uglich $x$.
Damit ist $\big(x_*(\cdot),u_*(\cdot)\big)$ ein starkes lokales Minimum. \hfill $\square$}
\end{beispiel}

\begin{beispiel}
{\textrm Im Beispiel \ref{BeipielRWUnendlich} mit Budgetbeschr"ankung,
\begin{eqnarray*}
&& J\big(x(\cdot),z(\cdot),u(\cdot)\big)=\int_0^\infty e^{-\varrho t}\big(1-u(t)\big)x(t) \, dt \to \sup, \\
&& \dot{x}(t)=u(t)x(t), \; x(0)=1,\qquad \dot{z}(t)=e^{-\varrho t}x(t), \; z(0)=0, \; \lim_{t \to \infty} z(t) = Z, \\
&& u \in [0,1], \qquad \varrho \in (0,1),
\end{eqnarray*}
ist jeder zul"assige Steuerungsprozess $\big(x(\cdot),z(\cdot),u(\cdot)\big)$ global optimal.
Die Voraussetzungen in Theorem \ref{SatzHBPMPUH} an einen zul"assigen Steuerungsprozess sind genau dann erf"ullt,
wenn die Steuerung $u(\cdot)$ dem Raum $L_1(\R_+,[0,1])$ angeh"ort.
Weiterhin sind f"ur jeden zul"assigen Steuerungprozess mit den Multiplikatoren
$$\lambda_0=1, \qquad p(t)= e^{-\varrho t}, \qquad q(t)= \varrho-1$$
die notwendigen Bedingungen (\ref{SatzPAUH1})--(\ref{SatzPAUH3}) erf"ullt.
Weiterhin ist die Hamilton-Funktion $\mathscr{H}$ offenbar konkav bez"uglich $(x,z)$ und es gilt Theorem \ref{SatzHBPMPUH}.} \hfill $\square$
\end{beispiel}

\begin{beispiel}[Abbau einer erneuerbaren Ressource] \label{ExampleErneuRessource}\index{Ressourcenabbau}
{\textrm Wir betrachten die Aufgabe
\begin{eqnarray}
&& \label{ErneuRessource1} J\big(x(\cdot),u(\cdot)\big)
   =\int_0^\infty e^{-\varrho t}\big[\pi\big(x(t)\big)-\kappa\big(x(t)\big)\big] u(t) \, dt \to \sup, \\
&& \label{ErneuRessource2} \dot{x}(t) = G\big(x(t)\big) - u(t), \quad x(0)=x_0>0,\quad u(t) \in [0, u_{\max}].
\end{eqnarray}
"Okonomische Interpretation: Durch die Funktion $G$ wird die nat"urliche dynamische Entwicklung der Ressource (oder einer Population) beschrieben,
wie diese sich ohne externe Einfl"usse entwickelt.
Diese nat"urliche Entwicklung wird durch die Abbaurate $u$ gest"ort.
Entsprechend des Angebotes $u$ wird ein Gewinn in H"ohe des Preises $\pi$ abz"uglich der Kosten $\kappa$ erzielt.
Dabei nehmen wir, dass Preis und Kosten steigen desto seltener die Ressource ist. \\[2mm]
In vielen Aufgaben dieser Form zeigt die L"osung ein station"ares Verhalten bzw. es tendiert die L"osung zu einem station"aren Paar $(\overline{x},\overline{u})$.
Eine solches station"ares Paar in Verbindung mit der Adjungierten $p(\cdot)$ nennen wir eine Gleichgewichtsl"osung\index{Gleichgewichtsl"osung}.
Wir identifizieren eine Gleichgewichtsl"osung $\big(\overline{x},\overline{u},p(\cdot)\big)$ durch folgende Eigenschaften:
$$\varphi(\overline{x},\overline{u}) = 0, \quad H_u\big(t,\overline{x},\overline{u},1,p(t)\big) =0, \quad
  \dot{p}(t) = -G'(\overline{x}) p(t) - e^{-\varrho t}\big(\pi'(\overline{x})-\kappa'(\overline{x})\big)\overline{u}.$$
Aus den ersten beiden Beziehung ergeben sich f"ur die Gleichgewichtsl"osung $\big(\overline{x},\overline{u},p(\cdot)\big)$:
$$G(\overline{x}) = \overline{u}, \qquad p(t)= e^{-\varrho t}\big(\pi(\overline{x})-\kappa(\overline{x})\big), \qquad \dot{p}(t)=-\varrho p(t).$$
Mit der Bezeichung $\gamma(x)= \pi(x)-\kappa(x)$ f"ur den Gewinn pro Einheit erhalten wir ferner in der adjungierten Gleichung:
\begin{eqnarray*}
\dot{p}(t)= -\varrho p(t) &=& -G'(\overline{x}) p(t) - e^{-\varrho t} \gamma'(\overline{x})\overline{u}, \\
-\varrho e^{-\varrho t} \gamma(\overline{x}) &=& -G'(\overline{x}) e^{-\varrho t} \gamma(\overline{x}) - e^{-\varrho t} \gamma'(\overline{x})G(\overline{x}), \\
0 &=& \big(\varrho -G'(\overline{x})\big) \gamma(\overline{x}) - \gamma'(\overline{x})G(\overline{x})
\end{eqnarray*}
oder ausgedr"uckt in der Form der Gewinnelastizit"at
$$\varepsilon_{\gamma,x}:=\overline{x} \frac{\gamma'(\overline{x})}{\gamma(\overline{x})} =\overline{x} \frac{\varrho -G'(\overline{x})}{G(\overline{x})}.$$
F"ur das logistische Wachstum $G(x)=x(r-Kx)$ ergibt sich die Gewinnelastitzit"at
$$\varepsilon_{\gamma,x} = -\frac{(r-K\overline{x})-(\varrho -K\overline{x})}{r-K\overline{x}} \in (-1,0)$$
im Fall einer initialen Wachstumsrate $G'(0)=r > \varrho$. \hfill $\square$}
\end{beispiel}

\begin{bemerkung}{\textrm
Die Darstellung der Pontrjagin-Funktion und der notwendigen Bedingungen im Maximumprinzip \ref{SatzPAUH} nennt man auch
Gegenwartswert-Schreibweise (``present value''). \index{Pontrjaginsches Maximumprinzip!oekonomische@--, Gegenwartswert-Schreibweise}
Durch die Einf"uhrung Momentanwert-Schreibweise (``current value'') \index{Pontrjaginsches Maximumprinzip!oekonomische@--, Momentanwert-Schreibweise}
mittels der Funktionen $q=e^{\varrho t}p$ und $\tilde{H}=e^{\varrho t}H$,
\begin{eqnarray*}
e^{\varrho t} H(t,x,u,p,1) &=& e^{\varrho t} \big[\langle p, \varphi(t,x,u) \rangle- e^{- \varrho t} f(t,x,u)\big] \\
                           &=& \langle q, \varphi(t,x,u) \rangle - f(t,x,u) = \tilde{H}(t,x,u,q,1),
\end{eqnarray*}
lassen sich im letzten Beispiel die Gleichgewichtsbedingungen in folgende Form "uberf"uhren:
$$\varphi(\overline{x},\overline{u}) = 0, \quad \tilde{H}_u(t,\overline{x},\overline{u},1,\overline{q}) =0, \quad
  \varrho \overline{q} - \tilde{H}_x(t,\overline{x},\overline{u},1,\overline{q})=0.$$
In der Momentanwert-Schreibweise kommt der Charakter eines Gleichgewichts vollst"andig zum Ausdruck.
Deswegen findet man in den Untersuchungen zu Gleichgewichten in der Literatur h"aufig das Maximumprinzip in der Momentanwert-Schreibweise. \hfill $\square$}
\end{bemerkung}

\begin{beispiel} \index{Kapitalakkumulation}
{\textrm Im Folgenden befassen wir uns mit einem neoklassischen Modell der "okonomischen Wachstumstheorie.
Diese Problemklasse geht auf die bereits erw"ahnte Arbeit von Ramsey \cite{Ramsey} zur"uck.
Es bezeichne $K(\cdot)$ das Verm"ogen einer "Okonomie,
$Y(\cdot)$ das Nationaleinkommen und $C(\cdot)$ die Konsumption.
Weiterhin wird zu jedem Zeitpunkt das Nationaleinkommen $Y$ in Konsumption und Investition aufgeteilt, d.\,h. $Y(t)=C(t)+\dot{K}(t)$.
Mit der Nutzenfunktion $U$ ergibt sich die allgemeine Aufgabenstellung
$$J\big(K(\cdot),C(\cdot)\big) = \int_0^\infty e^{-\varrho t} U\big(C(t)\big) \, dt \to \sup, \quad \dot{K}(t) = Y(t)-C(t).$$
Das Nationaleinkommen bestimmt sich durch die produzierten G"utern,
die aus dem eingesetzen Kapital $K$ und den Arbeitsressourcen $L$ gewonnen werden.
Mit der Cobb-Douglas-Produktionsfunktion $F(K,L)$ \index{Funktion, absolutstetige!Cobb@--, Cobb-Douglas-Produktions-} \index{Cobb-Douglas-Produktionsfunktion}
erhalten wir die Darstellung
$$Y(t)=F\big(K(t),L(t)\big), \qquad F(K,L)=K^\alpha L^{1-\alpha}, \qquad \alpha \in (0,1).$$
Die Verteilung des Nationaleinkommens $Y$ in Investition und Konsumption beschrieben mit $u \in [0,1]$ liefert
$$\dot{K}(t) = u(t) F\big(K(t),L(t)\big), \qquad C(t)=\big(1-u(t)\big) F\big(K(t),L(t)\big).$$
Mit einer logarithmischen Nutzenfunktion $U$ stellt sich so die Aufgabe wie folgt dar:
\begin{eqnarray*}
&& J\big(K(\cdot),C(\cdot)\big) = \int_0^\infty e^{-\varrho t} \ln\big[ \big(1-u(t)\big) F\big(K(t),L(t)\big) \big] \, dt \to \sup, \\
&& \dot{K}(t) = u(t) F\big(K(t),L(t)\big), \qquad \dot{L}(t) = \mu L(t), \qquad u(t) \in [0,1].
\end{eqnarray*}
Wir f"uhren die Kapital- und Konsumintensit"at $k = K/L$, $c=C/L$ und die Produktionsfunktion $f(x)=F(x,1)$ ein.
Dann ergeben sich die Beziehungen
\begin{eqnarray*}
\dot{k}(t) &=& \frac{\dot{K}(t) L(t)- K(t)\dot{L}(t)}{L^2(t)} = \frac{u(t) K^\alpha(t) L^{1-\alpha}(t) L(t)- K(t)\dot{L}(t)}{L^2(t)} \\
           &=& F\big(K(t)/L(t),1\big) - \frac{\dot{L}(t)}{L(t)}\frac{K(t)}{L(t)} =u(t)f\big(k(t)\big)-\mu k(t), \\
c(t)       &=& \big(1-u(t)\big) \frac{F\big(K(t),L(t)\big)}{L(t)} = \big(1-u(t)\big) F\big(K(t)/L(t),1\big) = \big(1-u(t)\big)f\big(k(t)\big).
\end{eqnarray*}
Auf diese Weise bekommt die Aufgabe die finale Form
\begin{eqnarray*}
&& J\big(k(\cdot),c(\cdot)\big) = \int_0^\infty e^{-\varrho t} \ln\big[ \big(1-u(t)\big) f\big(k(t)\big) \big] \, dt \to \sup, \\
&& \dot{k}(t) = u(t) f\big(k(t)\big)-\mu k(t), \qquad u(t) \in [0,1].
\end{eqnarray*}
Die Pontrjagin-Funktion der Aufgabe in der Momentanwert-Schreibweise ist
$$\tilde{H}(t,k,u,q,1) = q[u f(k)-\mu k] + \ln[(1-u) f(k)].$$
Dann ergeben die Gleichgewichtsbedingungen in der Momentanwert-Schreibweise:
\begin{enumerate}
\item[(a)] Das dynamische Gleichgewicht liefert
           $$\varphi(\overline{k},\overline{u}) = \overline{u} f(\overline{k})-\mu \overline{k} = 0
             \qquad\Leftrightarrow\qquad \overline{u}  = \mu\frac{\overline{k}}{f(\overline{k})}.$$ 
\item[(b)] Der stabile Verteilungsparameter erf"ullt die Beziehung
           $$\tilde{H}_u(t,\overline{k},\overline{u},1,\overline{q}) = \overline{q} f(\overline{k}) - \frac{1}{1-\overline{u}} = 0
             \qquad\Leftrightarrow\qquad \overline{u}  = 1- \frac{1}{\overline{q}f(\overline{k})}.$$
\item[(c)] F"ur den konstanten Schattenpreis erhalten wir in der adjungierten Gleichung
           $$\varrho \overline{q} - \tilde{H}_k(t,\overline{k},\overline{u},1,\overline{q})
             = \varrho \overline{q} -  \overline{q}[\overline{u} f(\overline{k})-\mu] - \frac{f'(\overline{k})}{f(\overline{k})} = 0.$$
           Die Division durch $\overline{q}$ und das Einsetzen der Beziehungen (a) und (b) liefert
           $$f'(\overline{k}) = \varrho+\mu.$$
\end{enumerate}
Wegen $f'(\overline{k}) > \mu$ ergibt sich $\overline{u} \in (0,1)$.
Damit wird durch die Beziehungen (a)--(c) ein eindeutiges Gleichgewicht $(\overline{k},\overline{u})$ festgelegt.
Abschlie"send bemerken wir,
dass die Hamilton-Funktion $\mathscr{H}$ stets konkav in $k$ ist. \hfill $\square$}
\end{beispiel}
       \subsection{Optimalit"atsbedingungen unter Zustandsbeschr"ankungen}
\begin{theorem}[Pontrjaginsches Maximumprinzip] \label{SatzPAUHZA}
Es sei $\big(x_*(\cdot),u_*(\cdot)\big) \in \mathscr{B}_{\textrm adm} \cap \mathscr{B}_{\textrm Lip}$.
Weiterhin nehmen wir an,
dass
$$\int_0^\infty \big\|\varphi\big(t,x_*(t),u_*(t)\big)\big\| \, dt < \infty, \qquad \int_0^\infty \big\|\varphi_x\big(t,x_*(t),u_*(t)\big)\big\| \, dt < \infty$$
ausfallen und es m"oge zu jedem $\delta>0$ ein $T>0$ existieren mit
\begin{eqnarray*}
&& \int_T^\infty \big\| \varphi\big(t,x(t),u_*(t)\big)-\varphi\big(t,x'(t),u_*(t)\big) - \varphi_x\big(t,x_*(t),u_*(t)\big)\big(x(t)-x'(t)\big) \big\| \, dt \\
&& \hspace*{20mm} \leq \delta \|x(\cdot)-x'(\cdot)\|_\infty
\end{eqnarray*}
f"ur alle $x(\cdot), x'(\cdot) \in W^1_\infty(\R_+,\R^n)$ mit $\|x(\cdot)-x_*(\cdot)\|_\infty < \gamma$, $\|x'(\cdot)-x_*(\cdot)\|_\infty < \gamma$. \\[2mm]
Ist $\big(x_*(\cdot),u_*(\cdot)\big)$ ein starkes lokales Minimum der Aufgabe (\ref{PAUH1})--(\ref{PAUH5}),
dann existieren eine Zahl $\lambda_0 \geq 0$, eine Vektorfunktion $p(\cdot):\R_+ \to \R^n$
und auf den Mengen
$$T_j=\big\{t \in \overline{\R}_+ \,\big|\, g_j\big(t,x_*(t)\big)=0\big\}, \quad j=1,...,l,$$
konzentrierte nichtnegative regul"are Borelsche Ma"se $\mu_j$ endlicher Totalvariation
(wobei s"amtliche Gr"o"sen nicht gleichzeitig verschwinden) derart, dass
\begin{enumerate}
\item[(a)] die Vektorfunktion $p(\cdot)$ von beschr"ankter Variation ist, der adjungierten Gleichung\index{adjungierte Gleichung}
           \begin{eqnarray}
           p(t)&=&- \lim_{t \to \infty}h_{1x}^T\big(t,x_*(t)\big) l_1 + \int_t^\infty H_x\big(s,x_*(s),u_*(s),p(s),\lambda_0\big) \, ds \nonumber \\
           \label{SatzPAUHZA1} & & -\sum_{j=1}^l \int_t^\infty g_{jx}\big(s,x_*(s)\big)\, d\mu_j(s)
           \end{eqnarray}
           gen"ugt und die Transversalit"atsbedingung\index{Transversalit"atsbedingungen}
           \begin{equation}\label{SatzPAUHZA2}
           p(0) = {h_0'}^T\big(x_*(0)\big)l_0
           \end{equation}
           erf"ullt;
\item[(b)] f"ur fast alle $t\in \R_+$ die Maximumbedingung
           \begin{equation}\label{SatzPAUHZA3}
           H\big(t,x_*(t),u_*(t),p(t),\lambda_0\big) = \max_{u \in U} H\big(t,x_*(t),u,p(t),\lambda_0\big)
           \end{equation}
           gilt.
\end{enumerate}
\end{theorem}

Die Adjungierte $p(\cdot)$ ist linksseitig stetig und es gilt
im rechten Endpunkt die Transversalit"atsbedingung 
$$\lim_{t \to \infty} p(t) = \lim_{t \to \infty} \bigg[-{h_1'}^T\big(x_*(t)\big)l_1-\sum_{j=1}^l g_{jx}\big(t,x_*(t)\big)\, \mu_j(\{\infty\})\bigg].$$

Es bezeichnet $\mathscr{H}(t,x,p) = \sup\limits_{u \in U} H(t,x,u,p,1)$ die Hamilton-Funktion $\mathscr{H}$ im normalen Fall.

\begin{theorem} \label{SatzHBPMPUHZB}
In der Aufgabe (\ref{PAUH1})--(\ref{PAUH5}) mit expliziten Randbedingungen sei
$\big(x_*(\cdot),u_*(\cdot)\big) \in \mathscr{B}_{\textrm Lip} \cap \mathscr{B}_{\textrm adm}$.
Au"serdem sei die Vektorfunktion $p(\cdot):\R_+ \to \R^n$ st"uckweise stetig,
besitze h"ochstens abz"ahlbar viele Sprungstellen $s_k \in (0,\infty)$,
die sich nirgends im Endlichen h"aufen,
und $p(\cdot)$ sei zwischen diesen Spr"ungen stetig differenzierbar. 
Ferne gelte:
\begin{enumerate}
\item[(a)] Das Tripel $\big(x_*(\cdot),u_*(\cdot),p(\cdot)\big)$
           erf"ullt (\ref{SatzPAUH1})--(\ref{SatzPAUH3}) in Theorem \ref{SatzPAUH} mit $\lambda_0=1$.        
\item[(b)] F"ur jedes $t \in \R_+$ ist die Funktion $\mathscr{H}\big(t,x,p(t)\big)$ konkav 
           und es sind die Funktionen $g_j(t,x)$, $j=1,...,l$, konvex bez"uglich $x$ auf $V_\gamma(t)$.
\end{enumerate}
Dann ist $\big(x_*(\cdot),u_*(\cdot)\big)$ ein starkes lokales Minimum der Aufgabe (\ref{PAUH1})--(\ref{PAUH5}) mit expliziten Randbedingungen.
\end{theorem}

\begin{beispiel} \label{ExampleRessource2}
{\textrm Wir betrachten erneut das Beispiel \ref{ExampleRessource}: \index{Ressourcenabbau}
\begin{eqnarray*}
&& J\big(x(\cdot),y(\cdot),u(\cdot)\big) =\int_0^\infty e^{-\varrho t}\big[pf\big(u(t)\big)-ry(t)-qu(t)\big] \, dt \to \sup, \\
&& \dot{x}(t) = -u(t),\quad \dot{y}(t)=cf\big(u(t)\big), \quad x(0)=x_0>0,\quad y(0)=y_0\geq 0, \\
&& x(t) \geq 0, \qquad u(t) \geq 0, \qquad b,c,\varrho, q,r >0, \qquad \varrho - rc>0.
\end{eqnarray*}
Wegen der Zustandsbeschr"ankung sind f"ur jeden zul"assigen Steuerungprozess die Beschr"ankungen an die Dynamik in Theorem \ref{SatzPAUHZA} erf"ullt.
Bei der Anwendung der notwendigen Bedingungen ist nun zu beachten,
dass das Gewicht $\nu(t)=e^{-at}$ in die Darstellung der Adjungierten nicht einflie"st,
d.\,h. mit $\lambda_0=1$ gen"ugt die Adjungierte $p_1(\cdot)$ der Gleichung
$$p_1(t)=\int_t^\infty \, d\mu(s).$$
Die weitere Diskussion kann nun vollst"andig aus Beispiel \ref{ExampleRessource} "ubernommen werden.
Insbesondere enth"alt die Pontrjagin-Funktion
$$H(t,x,y,u,p_1,p_2,1) = p_1 (-u)+p_2 cf(u) + e^{-\varrho t}[f(u)-ry-qu]$$
bez"uglich der Zustandsvariablen $x$ und $y$ nur den Term $e^{-\varrho t}ry$.
Damit ist die Hamilton-Funktion offenbar konkav bez"uglich $(x,y)$.
Au"serdem ist die Zustandsbeschr"ankung linear in $x$.} \hfill $\square$
\end{beispiel}

\addcontentsline{toc}{section}{Literatur}
\lhead[\thepage \, Literatur]{Optimale Steuerung mit unendlichem Zeithorizont}
\rhead[Optimale Steuerung mit unendlichem Zeithorizont]{Literatur \thepage}

\end{document}